\font\twelvemsb=msbm10 at 12pt
\font\eightmsb=msbm8
\newfam\msbfam
\textfont\msbfam=\twelvemsb
\scriptfont\msbfam=\eightmsb
\def\Bbb#1{\fam\msbfam\relax#1}

\def\endproof{\ \  \rule{0.5em}{0.5em}}

\documentstyle[12pt]{article}
\advance\textheight by \topskip
\newtheorem{theorem}{Theorem}[section]
\newtheorem{proposition}[theorem]{Proposition}

\newcommand{\N}{{\Bbb{N}}}
\newcommand{\R}{{\Bbb{R}}}

\newcommand{\Rpp}{{\Bbb{R}}^{p \times p}}
\newcommand{\Rpq}{{\Bbb{R}}^{p \times q}}

\begin{document}
\title{Polynomial interpolation and Gaussian quadrature for matrix valued
functions}
\author{Ann Sinap\thanks{e-mail fgaee04@cc1.KULeuven.ac.be}\quad  and \quad
        Walter Van Assche\thanks{e-mail fgaee03@cc1.KULeuven.ac.be;
        Research Asscociate of the Belgian National Fund for Scientific
        Research} \\
        Department of Mathematics \\
        Katholieke Universiteit Leuven \\
        Celestijnenlaan 200\thinspace B \\
        B--3001 Heverlee, BELGIUM}
\date{}
\maketitle
\begin{abstract}
The techniques for polynomial interpolation and Gaussian
quadrature are generalized to matrix-valued functions. It is shown how the
zeros and rootvectors of matrix orthonormal polynomials can be used to get a
quadrature formula with the highest degree of precision.
\end{abstract}

\section{Introduction}
The aim of this paper is to construct quadrature formulas, using orthogonal
matrix polynomials, to approximate matrix integrals.
We will give an expression for the quadrature coefficients and show that
the constructed formula has the highest possible degree of accuracy and
converges to the exact value of the matrix integral.
All these ideas are generalizations of the classical Gaussian quadrature
rules for the scalar case.

In Section 2 we will give a survey
of definitions and properties of matrix polynomials. These can be found in the
book on matrix polynomials by I.\ Gohberg, P.\ Lancaster and L.\ Rodman
\cite{go:1982} and in the survey on orthogonal matrix polynomials by L.\
Rodman \cite{ro:1990}. In Section 3 we introduce orthogonal matrix
polynomials on the real line and discuss some properties which we will need in
the next sections. These orthogonal matrix polynomials have been considered
earlier by Delsarte, Genin and
Kamp \cite{de:1978} and Geronimo \cite{ge:1981}, \cite{ge:1982}.
As in the scalar case, the theory of approximate integration uses results from
the theory of interpolation. In Section 4 we will discuss polynomial
interpolation and particularly the interpolation problem of Lagrange. We
will give an expression for the Lagrange interpolation polynomial
in the general case and then apply this result to the case of a Jordan
pair $(X,J)$ of the orthonormal matrix polynomial $P_n(x)$. These results
are also generalizations of the known results for the scalar case. The
interpolation problem has also been treated in \cite{goka:1981},
\cite{goka:1982}, \cite{fu:1983} and \cite{baka:1990}, but their approach is an
algebraic one. We have restricted our attention to the interpolation
formula we needed for the construction of Gaussian quadrature rules.
The Gaussian quadrature formula is then constructed in Section 5, where we give
a formula for the quadrature coefficients and show that the quadrature
rule converges under appropriate conditions.

\section{Matrix polynomials}
If $A_{0},A_{1},\ldots,A_{n}$ are elements of $\Rpp$ and $A_{n} \not = 0 $, then
we call
$$ P(x) = A_{n}x^{n} + A_{n-1}x^{n-1} + \ldots + A_{1}x + A_{0}, $$
a {\it matrix polynomial} of degree $n$. This matrix polynomial is monic
when $A_{n} =I$, the identity matrix.
A point $x_{0}$ is a {\it zero} of $P(x)$, if $\det{P(x_{0})} = 0$. Note that if
the
leading coefficient of $P(x)$ is non-singular, $\det{P(x)}$ is a polynomial of
degree $np$.
Another important notion associated with matrix polynomials are Jordan chains. A
sequence of $p$-dimensional column vectors $v_{0},v_{1},\ldots,v_{k}$ is called
a {\it right Jordan chain} of length $k+1$ of a monic matrix polynomial
$\hat{P}(x)$ corresponding to $x_{0}$, if $v_{0} \not = 0 $ and
$$ \sum_{i=0}^{l} \frac{1}{i!} {\hat{P}}^{(i)}(x_{0}) v_{l-i} = 0, \qquad
l=0,1,\ldots,k.$$
The initial vector $v_{0} \not = 0 $ is called a
{\it rootvector} of $\hat{P}(x)$ corresponding to $x_{0}$.
Note that in \cite{go:1982} and \cite{ro:1990} the zeros are called
eigenvalues and the rootvectors are called eigenvectors.
${\hat{P}}^{(i)}(x)$ is  the $i$th derivative of $\hat{P}(x)$ with respect to
$x$ and this means that we take the $i$th derivative of every element of
$\hat{P}(x)$ with respect to $x$.

In analogy with the definition of a right Jordan chain, we call a
sequence of $p$-dimensional row vectors $w_{0},w_{1},\ldots,w_{k}$
a left Jordan chain of length $k+1$ of a monic matrix polynomial $\hat{P}(x)$
corresponding to $x_{0}$, if $v_{0} \not = 0 $ and
$$ \sum_{i=0}^{l} \frac{1}{i!} w_{l-i} {\hat{P}}^{(i)}(x_{0})  = 0, \qquad
l=0,1,\ldots,k.$$
Jordan chains are not unique: a matrix polynomial can have different
Jordan chains of various lengths and different rootvectors,
corresponding to the same zero. In what follows we formulate the
definitions and properties for right Jordan chains, unless explicitly mentioned.

A convenient way of writing a Jordan chain is given by the following property.
\begin{proposition}[\cite{go:1982}, p.\ 27]
The vectors $v_{0},v_{1},\ldots,v_{k}$ form a right Jordan chain of the
monic matrix polynomial
$ \hat{P}(x) = Ix^{n} + A_{n-1}x^{n-1} + \ldots + A_{1}x + A_{0} $
corresponding to $x_{0}$ if and only if $v_{0} \not = 0$ and
$$ X_{0}J_{0}^{n} + A_{n-1}X_{0}J_{0}^{n-1} + \ldots + A_{1}X_{0}J_{0} +
A_{0}X_{0} = 0, $$
where $X_{0} = \left( \begin{array}{ccc}
   v_{0}  &  \ldots &   v_{k}   \\
                      \end{array} \right)   $ is  a $p \times (k+1)$ matrix and
$J_{0}$ is a Jordan block of size $ (k+1) \times (k+1)$ with $x_{0}$ on the
main diagonal.
\end{proposition}

Observe that the equations
$$ \sum_{i=0}^{l} \frac{1}{i!} {L}^{(i)}(x_{0}) v_{l-i} = 0, \qquad
l=0,1,\ldots,k,$$
where $L(x)=L_{n}x^{n}+L_{n-1}x^{n-1}+\ldots+L_{1}x+L_{0}$, with $ L_{i} \in
\Rpq$ (i=1,2,\ldots,n), can always be written as
$$L_{n}X_{0}J_{0}^{n} + L_{n-1}X_{0}J_{0}^{n-1} + \ldots + L_{1}X_{0}J_{0} +
L_{0}X_{0} = 0,  $$
where $X_{0} = \left( \begin{array}{ccc}
   v_{0}  &  \ldots &   v_{k}   \\
                      \end{array} \right)   $ is  a $p \times (k+1)$ matrix and
$J_{0}$ is a Jordan block of size $ (k+1) \times (k+1)$ with $x_{0}$ on the
main diagonal.

In the following definitions and properties we restrict ourselves to monic $p
\times p$ matrix polynomials, but most of the theory can also be given in the
context of regular matrix polynomials. These are matrix polynomials which
satisfy $\det{P(x)} \not \equiv 0$.

Now we introduce the notion of a {\it canonical set of Jordan chains}.
Let
$$ v_{j,0}^{(i)}, v_{j,1}^{(i)},\ldots,v_{j,\mu_{j}^{(i)}-1}^{(i)}, \qquad
j=1,2,\ldots,s_{i},$$ be a set of Jordan chains of a $p \times p$
monic
matrix polynomial $\hat{P}(x)$ corresponding to the zero $x_{i}$. Then we
call the set canonical if the rootvectors $ v_{1,0}^{(i)}, v_{2,0}^{(i)},
\ldots, v_{s_{i},0}^{(i)}$ are linearly independent and $\sum_{j=1}^{s_{i}}
\mu_{j}^{(i)} = m_{i}$, where $m_{i}$ is the multiplicity of $x_{i}$ as zero of
$\hat{P}(x)$.
Such a canonical set of Jordan chains is not unique, but the number of
chains and their length depend only upon $\hat{P}(x)$ and $x_{i}$  and do not
depend on the choice of canonical set. A canonical set of Jordan chains
can be associated with a pair of matrices  $(X_{i},J_{i})$, which is called the
Jordan pair of $\hat{P}(x)$ corresponding to $x_{i}$ and defined as follows :
$$  X_{i} = ( v_{1,0}^{(i)}\ \ldots \ v_{1,\mu_{1}^{(i)}-1}^{(i)} \ \ldots \
 v_{s_{i},0}^{(i)}\ \ldots \ v_{s_{i},\mu_{s_{i}}^{(i)}-1}^{(i)}) \quad {\rm a\
} p \times m_{i} \ {\rm dimensional\ matrix}$$ and
$$ J_{i} = diag(J_{i,1},J_{i,2},\ldots,J_{i,s_{i}}) \quad {\rm a\ } m_{i}
\times m_{i}\  {\rm dimensional\ matrix\ with\ }$$
$$ J_{i,j} = \left ( \begin{array}{ccccc}
         x_{i} &   1    &        &        &       \\
               & \ddots & \ddots &        &       \\
               &        & \ddots & \ddots &       \\
               &        &        & x_{i}  &   1    \\
               &        &        &        &  x_{i}
                \end{array} \right ) \quad
{\rm a\ }  \mu_{j}^{(i)} \times  \mu_{j}^{(i)} \  {\rm dimensional\ matrix.} $$
A pair of matrices $(X,J)$ where $X$ is a $p \times np$ dimensional matrix  and
$J$ a $np \times np $ dimensional Jordan matrix is called a {\it Jordan
pair} for the monic matrix polynomial $\hat{P}(x)$ if
$$ X = ( X_{1}\ X_{2} \ \ldots \ X_{k} ) \quad {\rm and } \quad
J = diag( J_{1}\ J_{2} \ \ldots \ J_{k} ), $$
where $(X_{i},J_{i})$ is a Jordan pair of $\hat{P}(x)$ corresponding to $x_{i}$
and  $k$ is the number of different zeros of $\hat{P}(x)$. Jordan pairs have the
following important property :
\begin{proposition}[\cite{go:1982}, p.\ 45]
Let $(X,J)$ be a pair of matrices where $X$ is of size $p \times np$ and $J$
is a Jordan matrix of size $np \times np $. Then $(X,J)$ is a Jordan pair of the
monic matrix polynomial
 $ \hat{P}(x) = Ix^{n} + A_{n-1}x^{n-1} + \ldots + A_{1}x + A_{0} $
if and only if
\begin{itemize}
\item[(1)] $${col(XJ^{l})}_{l=0}^{n-1} =
\left( \begin{array}{c}
       X \\ XJ \\ XJ^{2} \\ \ldots \\ XJ^{n-1}
   \end{array} \right )  $$
is a non-singular $np \times np$ matrix,
\item[(2)] $ XJ^{n} + A_{n-1}XJ^{n-1} + \ldots + A_{1}XJ + A_{0}X = 0 $.
\end{itemize}
\end{proposition}

The requirement that $J$ is Jordan is not essential. We call a pair of matrices
$(X,T)$ where $X$ is of size $p \times np$ and $T$
is a $np \times np $ dimensional matrix, a {\it standard pair} for the
monic matrix
polynomial $ \hat{P}(x) = Ix^{n} + A_{n-1}x^{n-1} + \ldots + A_{1}x + A_{0} $
if
\begin{itemize}
\item[(1)] ${col(XT^{l})}_{l=0}^{n-1}$ is a non-singular $np \times np$ matrix,
\item[(2)] $ XT^{n} + A_{n-1}XT^{n-1} + \ldots + A_{1}XT + A_{0}X = 0 $.
\end{itemize}
This means that every Jordan pair is a standard pair and every standard pair
$(X,T)$ for which $T$ is a Jordan matrix is a Jordan pair.

With every standard pair $(X,T)$ of a monic matrix polynomial $\hat{P}(x)$ we
can associate a third matrix $Y$ of size $np \times p$, with
$$ Y = {\left( \begin{array}{c}
       X \\ XJ \\ XJ^{2} \\ \vdots \\ XJ^{n-1}
   \end{array} \right )}^{-1} \,
        \left( \begin{array}{c}
       0 \\ 0  \\ \vdots \\ 0  \\ I
   \end{array} \right ). $$
The triple $(X,T,Y)$ is called a {\it standard triple} for $\hat{P}(x)$
and if $T=J$,
a Jordan matrix, then $(X,J,Y)$ is called a {\it Jordan triple}. One can
proof that $ ( X^{\prime} , C_{1} , Y^{\prime} )$, where
$$ X^{\prime} = \left( \begin{array}{cccc}
I &  0  &  \ldots & 0   \end{array}  \right), \ \
C_{1} = \left( \begin{array}{ccccc}
   0  &  I  &  0  &  \ldots  &  0  \\
   0  &  0  &  I  &  \ldots  &  0  \\
\vdots &  \vdots &  \vdots  &  \ddots  &  0  \\
   0  &  0  &  0  & \ldots & I \\
  -A_{0} & -A_{1} & -A_{2} & \ldots & -A_{n-1} \\
  \end{array} \right) \quad {\rm and } \ \ \
Y^{\prime}   =  \left( \begin{array}{c}
       0 \\ 0  \\ \ldots \\ 0  \\ I
   \end{array} \right ), $$
is a standard triple of
$\hat{P}(x)= Ix^{n} + A_{n-1}x^{n-1} + \ldots + A_{1}x + A_{0} $.
Two standard triples $(X,T,Y)$ and
$ ( X^{\prime},T^{\prime}, Y^{\prime} )$ are {\it similar} if there exists
an invertible $np \times np $ matrix $S$ such that
$$ X^{\prime}=XS \quad ; \quad
T^{\prime}=S^{-1}TS \quad {\rm and} \quad Y^{\prime}=S^{-1}Y. $$
This matrix $S$ is uniquely defined by,
$$ S = {({col(XT^{l})}_{l=0}^{n-1})}^{-1} \,
{col(X^{\prime}{T^{\prime}}^{l})}_{l=0}^{n-1}. $$

Consider a Jordan triple $(X,J,Y)$, then we already know that the columns
of $X$, when decomposed into blocks consistently with the decomposition of $J$
into blocks, form right Jordan chains for $\hat{P}(x)$. In analogy with this, we
can give a similar meaning to the rows of $Y$. Indeed, the rows of $Y$,
partitioned into blocks consistently with the decomposition of $J$
into blocks and taken in reverse order, form left Jordan chains for
$\hat{P}(x)$.
The notion of standard triples is important for the following representation
theorem:
\begin{theorem}[\cite{go:1982}, p.\ 58]
Let $ \hat{P}(x) = Ix^{n} + A_{n-1}x^{n-1} + \ldots + A_{1}x + A_{0} $
be a monic matrix polynomial of degree $n$ with standard triple $(X,T,Y)$. Then
$\hat{P}(x)$ admits the following representations :
\begin{itemize}
\item[(1)] $\hat{P}(x) = x^{n}I - XT^{n} (V_{1}+V_{2}x+\ldots+V_{n}x^{n-1}),
$ where $V_{i}$ are $np \times p $ matrices such that
$ (V_{1} \ V_{2}\ \ldots \ V_{n} ) =
{ ( {col(XT^{l})}_{l=0}^{n-1} )}^{-1} $.
\item[(2)]
$\hat{P}(x) = x^{n}I - (W_{1}+W_{2}x+\ldots+W_{n}x^{n-1})T^{n}Y, $  where
$W_{i}$ are $p \times np $ matrices such that
$ {col(W_{l})}_{l=1}^{n} = {\left( \begin{array}{ccccc}
Y &  TY & T^{2}Y & \ldots & T^{n-1}Y
  \end{array} \right) }^{-1} $. \end{itemize}
\end{theorem}
Note that those forms are independent of the choice of the standard triple.

Finally we have some properties about the divisibility of matrix polynomials.
We say that the matrix polynomials $Q(x)$ and $R(x)$ are the {\it right
quotient} and the {\it right remainder},
respectively, of $P(x)$ on division by $D(x)$ if $$ P(x) = Q(x) D(x) + R(x) $$
and if the degree of $R(x)$ is less than that of $D(x)$. The right division of
matrix polynomials of the same order is always possible and unique, provided the
divisor is a polynomial with non-singular leading coefficient (see \cite[p.
78]{ga:1960}).
In a similar manner we can define the left division : $ P(x) = D(x)
\tilde{Q}(x) + \tilde{R}(x)  $.

\begin{proposition}
Let $P(x)=A_{m}x^{m}+A_{m-1}x^{m-1}+\ldots+A_{1}x+A_{0}$ be a matrix polynomial
of degree $m$ and  let
 $D(x)=B_{n}x^{n}+B_{n-1}x^{n-1}+\ldots+B_{1}x+B_{0}$ be a matrix polynomial
 of degree $n$, with a non-singular leading coefficient $B_{n}$ and with Jordan
pair $(X,J)$. Then $D(x)$  is a right divisor of $P(x)$ if and only if
$$ A_{m}XJ^{m} + A_{m-1}XJ^{m-1} + \ldots + A_{1}XJ + A_{0}X = 0. $$
\end{proposition}

\noindent
\underline{Proof} \\
First of all we denote that every Jordan chain $v_{0},v_{1},\ldots,v_{k}$ of
$D(x)$ and corresponding to $x_{0}$ is also a Jordan chain of $Q(x)D(x)$,
corresponding to the same zero :
\begin{eqnarray*}
\sum_{i=0}^{l} \frac{1}{i!} {[Q(x)D(x)]}^{(i)}_{|x=x_{0}} v_{l-i}  & = &
\sum_{i=0}^{l} \sum_{t=0}^{i} \frac{1}{i!} \left( \begin{array}{c} i \\ t
\end{array} \right) Q^{(t)}(x_{0}) D^{(i-t)}(x_{0}) v_{l-i} \\
& = &
\sum_{t=0}^{l} \frac{1}{t!} Q^{(t)}(x_{0}) \sum_{i=t}^{l} \frac{1}{(i-t)!}
D^{(i-t)}(x_{0}) v_{l-i} \\
& = &
\sum_{t=0}^{l} \frac{1}{t!} Q^{(t)}(x_{0}) \sum_{s=0}^{l-t} \frac{1}{s!}
D^{(s)}(x_{0}) v_{l-t-s} \\
& = &   0 \quad {\rm for}  \quad l=0,1,\ldots,k.
\end{eqnarray*}
Suppose $D(x)$ is a right divisor of $P(x)$, this means $P(x)=Q(x)D(x)$. But
then every Jordan chain of $D(x)$ is also a Jordan chain of $P(x)$.
In particular we have
$$ A_{m}XJ^{m} + A_{m-1}XJ^{m-1} + \ldots + A_{1}XJ + A_{0}X = 0. $$
Suppose now we know that the above mentioned equation holds.
This means that every Jordan chain of the Jordan pair of $D(x)$ is also a
Jordan chain of $P(x)$. Moreover, since the leading coefficient of $D(x)$ is
non-singular, the right division is possible and unique, $P(x)=Q(x)D(x)+R(x)$.
So, every Jordan chain of the Jordan pair of $D(x)$ is also a Jordan chain of
$R(x)$. But this implies that the matrix polynomial $R(x)$ of degree $
\leq n-1$ satisfies $$ \left( \begin{array}{cccc}
R_{0} & R_{1} & \ldots & R_{n-1} \\
\end{array} \right) \,
\left( \begin{array}{c}
X  \\   XJ  \\  \vdots \\  XJ^{n-1}  \\
\end{array} \right) = 0. $$
Since $(X,J)$ is a Jordan pair, the $np \times np$ dimensional matrix
${col(XJ^{l})}_{l=0}^{n-1}$ is non-singular and thus $R(x) = 0$.
\endproof

Note that this proposition was also proved in \cite[Thm.\ 2.1 on p.\
333]{goka:1981}.

\section{ Orthogonal matrix polynomials on the real line}
A symmetric $p \times p$ matrix-valued function $W(x)$, integrable over
$[a,b]$
is called a weight matrix function if $W(x) \geq 0$ and $\det{W(x)} \not = 0$
holds almost everywhere (see \cite{de:1978}).
The notation $X \leq Y $ for symmetric matrices means that  $Y-X$ is positive
semidefinite.
The matrix integral
$$ \int_{a}^{b} \ F(x) \, W(x) \, G(x)^{T}\, dx $$
where $F(x)$ and $G(x)$ are continuous matrix-valued functions, is defined in a
natural way. The $(i,j)$th element is given by a sum of integrals :
$$ \sum_{s=1}^{p} \sum_{t=1}^{p} \int_{a}^{b} {F(x)}_{i,s}\, {W(x)}_{s,t}\,
{G(x)^{T}}_{t,j}\, dx. $$
Let $\Rpp [x]$ be the set of polynomials in a real variable $x$ and whose
coefficients are $p \times p$ matrices with real entries. If $P(x)$ and $Q(x)$
are elements of this set, then we define a matrical inner product on $\Rpp [x]$
as follows:
$$ {\langle P(x),Q(x) \rangle}_{L}\ =\ \int_{a}^{b} \ P(x) \, W(x) \,
 Q(x)^{T}\, dx. $$
This matricial inner product has some properties which we will recognize as
generalizations of the properties of the scalar inner product.
\begin{proposition} \label{eig:emi}
\begin{enumerate}
\item  ${\langle P,Q \rangle}_{L} \ = \
 {\langle Q,P \rangle_{L} }^{T} $ where $P,Q \in  \Rpp [x]$.
\item  $ {\langle C_{1}P_{1}+C_{2}P_{2},Q \rangle}_{L}  =
C_{1} {\langle P_{1},Q \rangle}_{L}+ C_{2} {\langle P_{2},Q \rangle}_{L} $
where $C_{1},C_{2} \in  \Rpp$ \ and $P_{1},P_{2},Q \in  \Rpp[x]$.
\item  ${\langle xP,Q \rangle}_{L} \ = \
 {\langle P,xQ \rangle_{L} } $ where $P,Q \in  \Rpp [x]$.
\item Let $P(x) \in \Rpp [x] $, then $ {\langle P,P \rangle }_{L}$ is
positive semidefinite and even
positive definite if $\det P(x) \not\equiv 0$.
\item  Let $P \in \Rpp [x]$, then
 ${ \langle P,P \rangle}_{L}=0 $ if and only if $P=0$.
\end{enumerate}
\end{proposition}
These properties are easily proved by means of straightforward computation.

A generalization of the Gram-Schmidt orthonormalisation procedure for the set \\
$ \{I,xI,x^{2}I,\ldots \}$ with respect to the matricial inner product
${\langle .,. \rangle}_{L}$ will give a set of orthonormal matrix polynomials
$\{ P_{n}(x) \}_{n=0}^{\infty}$ which satisfy
$$ \int_{a}^{b} P_{n}(x)\, W(x)\, {P_{m}(x)}^{T}\, dx  = \delta_{n,m}I. $$
Moreover, $P_{n}(x)$ is a matrix polynomial of degree $n$, with a
non-singular
leading coefficient and is defined upon a multiplication on the left by
an orthogonal matrix.

As in the scalar case, these orthonormal matrix polynomials are
orthogonal to every matrix polynomial of lower degree and they satisfy a
three-term recurrence relation.
$$ x P_{n}(x)=D_{n+1}P_{n+1}(x)+E_{n}P_{n}(x)+D_{n}^T P_{n-1}(x),\qquad
n \geq 0, $$
$$ P_{-1}(x)=0 \quad {\rm and} \quad P_{0}(x)=I, $$
where $D_{n}$ is a positive definite matrix and
$E_{n}$ is a  symmetric matrix.
The orthonormal polynomials are defined only up to a left orthogonal
factor and it is convenient to choose this factor in such a way
that the recurrence coefficients $D_n$ are symmetric.
We assumed, without loss of generality  that
$\int_{a}^{b} W(x) \, dx = I $. Furthermore we have the
Christoffel-Darboux formula :
$$ P_{n}(y)^{T} D_{n+1} P_{n+1}(x) - P_{n+1}(y)^{T} D_{n+1} P_{n}(x) =
 (x-y) \ \sum_{i=0}^{n} P_{i}(y)^{T} P_{i}(x).  $$
 (see \cite{ge:1982}).
If we take $x=y$ , we get
$$ P_{n}(x)^{T} D_{n+1} P_{n+1}(x) = P_{n+1}(x)^{T} D_{n+1} P_{n}(x)$$
so that $ P_{n}(x)^{T} D_{n+1} P_{n+1}(x) $
is a symmetric matrix.  By means of
straightforward computation we get the following equation :
 $$ \sum_{i=0}^{n} P_{i}(x)^{T} P_{i}(x)  =
{P_{n+1}^{T}(x)}^{\prime} D_{n+1} P_{n}(x) - {P_{n}^{T}(x)}^{\prime}
D_{n+1} P_{n+1}(x). $$
The matrix
$$ K_{n}(x,y) = \sum_{i=0}^{n} P_{i}(y)^{T} P_{i}(x)  $$ is a positive definite
matrix and we call it the reproducing kernel because of the following property.
\begin{proposition}
Let $ \Pi_{m}(x)$ be a matrix polynomial of degree  $m \leq n$, then
 $${\langle \Pi_{m}(x), K_{n}(x,y) \rangle}_{L} = {\Pi_{m}(y)} $$
\end{proposition}

\noindent
\underline{Proof} \\
If we write $ \Pi_{m}(x)$ in terms of the orthonormal matrix polynomials
$P_{0}(x),\ldots ,P_{m}(x)$ :  $$  \Pi_{m}(x) = \sum_{i=0}^{m} A_{i} \,
P_{i}(x), $$ we have for $m \leq n$
\begin{eqnarray*}
{\langle \Pi_{m}(x), K_{n}(x,y) \rangle}_{L} & =  &
  \sum_{i=0}^{m} \,  \sum_{j=0}^{n} A_{i}  \,
{\langle P_{i}(x),P_{j}(x) \rangle}_{L}\, P_{j}(y)  \\
& = &   \sum_{i=0}^{m}  A_{i}\, P_{i}(y)  \\
& = & {\Pi_{m}(y)}. \ \ \rule{0.5em}{0.5em}
\end{eqnarray*}
\bigskip

In the scalar case all the zeros of an orthonormal polynomial are simple.
This is not the case for orthonormal matrix polynomials, but nevertheless we can
proof a similar property.
\begin{proposition}
The zeros of the orthonormal matrix polynomial $P_{n}(x)$ have a
multiplicity $ \leq p$, where p is the size of the matrices.
\end{proposition}

\noindent
\underline{Proof} \\
Let $x_{0}$ be a zero of $P_{n}(x)$ with multiplicity $m > p$. Consider
a canonical set of right Jordan chains corresponding to $x_{0}$
$$ v_{i,0}, v_{i,1}, \ldots , v_{i,\mu_{i}-1}, \qquad  i=1,2,\ldots,s.
$$ This means that the $p$-dimensional column vectors
$ v_{1,0},v_{2,0},\ldots,v_{s,0} $ are linearly independent and
$\sum_{i=1}^{s} \mu_{i} = m$ (Section 1).
Since $m > p$, there has to be a Jordan chain of length $ > 1$. Suppose
$v_{0}$ and $v_{1}$ are the two leading vectors of this chain, then they satisfy
\begin{eqnarray*}
& & P_{n}(x_{0}) v_{0} = 0  \qquad  {\rm and} \qquad  v_{0} \not=0,  \\
& & P_{n}^{\prime}(x_{0}) v_{0} + P_{n}(x_{0}) v_{1} = 0.
\end{eqnarray*}
Using these equations, we get
\begin{eqnarray*}
v_{0}^{T} K_{n-1}(x_{0},x_{0}) v_{0}
& = & v_{0}^{T} {P_{n}^{T}(x_{0})}^{\prime} D_{n} P_{n-1}(x_{0}) v_{0} -
 v_{0}^{T} {P_{n-1}^{T}(x_{0})}^{\prime} D_{n} P_{n}(x_{0}) v_{0}   \\
& =  & - v_{1}^{T} P_{n}^{T}(x_{0}) D_{n} P_{n-1}(x_{0}) v_{0} \\
& =  & - v_{1}^{T} P_{n-1}^{T}(x_{0}) D_{n} P_{n}(x_{0}) v_{0} \\
& = & 0.
\end{eqnarray*}
But $K_{n-1}(x_{0},x_{0})$ is a symmetric and positive definite matrix and
$v_{0} \not= 0$. So the multiplicity of $x_{0}$ as zero of $P_{n}(x)$ has to
be $ \leq p$.
\endproof  \bigskip

\noindent
\underline{Corollary.} \\
In this proof we showed that the length of a Jordan chain of
$P_{n}(x)$ cannot be greater than $1$. Thus a canonical set of right Jordan
chains of  $P_{n}(x)$ corresponding with a zero $x_{0}$ consist of $m$
linearly independent, non-zero $p$-dimensional column vectors, where $m$ is the
multiplicity of $x_{0}$ as zero of $P_{n}(x)$.

\section{Polynomial interpolation}

\subsection{Polynomial interpolation in general}
Consider a $p \times p$ matrix-valued function $F(x)$ and $k$ different points
 $x_{1},x_{2},\ldots,x_{k}$ with multiplicity resp.\
$m_{1},m_{2},\ldots,m_{k}$ where $\sum_{i=1}^{k} m_{i} = np$. With every point
$x_{i}$, $i =1,2,\ldots,k$, we associate a set of $p$-dimensional column
vectors $$ v_{1,0}^{(i)}, v_{1,1}^{(i)},\ldots,v_{1,\mu_{1}^{(i)}-1}^{(i)},
 v_{2,0}^{(i)}, v_{2,1}^{(i)},\ldots,v_{2,\mu_{2}^{(i)}-1}^{(i)},\ldots,
 v_{s_{i},0}^{(i)}, v_{s_{i},1}^{(i)},\ldots,v_{s_{i},\mu_{s_{i}}^{(i)}-1}^{(i)}
$$
where $\sum_{j=1}^{s_{i}} \mu_{j}^{(i)} = m_{i} $ and
$ v_{1,0}^{(i)}, v_{2,0}^{(i)},\ldots,v_{s_{i},0}^{(i)} $ are non-zero, linearly
independent vectors. When we put these vectors
in a $p \times np$ dimensional matrix, we get
$$ X = ( X_{1}\ X_{2} \ \ldots \ X_{k} ) \quad {\rm where } \quad
   X_{i} = ( v_{1,0}^{(i)}\ \ldots \ v_{1,\mu_{1}^{(i)}-1}^{(i)} \ \ldots \
 v_{s_{i},0}^{(i)}\ \ldots \ v_{s_{i},\mu_{s_{i}}^{(i)}-1}^{(i)}). $$
The square $np \times np $ matrix $J$ is given by
$$ J = diag( J_{1},J_{2},\ldots, J_{k} ) \quad {\rm where } \quad
 J_{i} = diag( J_{i,1},J_{i,2},\ldots,J_{i,s_{i}}) $$
$$  \quad {\rm and } \quad
 J_{i,j} = \left ( \begin{array}{ccccc}
         x_{i} &   1    &        &        &       \\
               & \ddots & \ddots &        &       \\
               &        & \ddots & \ddots &       \\
               &        &        & x_{i}  &   1    \\
               &        &        &        &  x_{i}
                \end{array} \right ). $$
Here $J_{i}$ is a $m_{i} \times m_{i} $ matrix and
$J_{i,j}$ a $\mu_{j}^{(i)} \times \mu_{j}^{(i)} $ matrix.
We require that the matrix
$ {col(XJ^{l})}_{l=0}^{n-1} $
is non-singular. Note that this is a necessary and sufficient condition in order
that $(X,J)$ is a Jordan pair of a monic matrix polynomial of degree $n$
(Section 1).

\noindent
We will construct a matrix polynomial $P(x)=A_{n-1}x^{n-1}+
\ldots+A_{1}x+A_{0} $, which satisfies
$$ \sum_{l=0}^{q} \frac{1}{l!} \, P^{(l)}(x_{i}) \, v_{t,q-l}^{(i)} \ = \
 \sum_{l=0}^{q} \frac{1}{l!} \, F^{(l)}(x_{i}) \, v_{t,q-l}^{(i)}  $$
and this for $$ q=0,1,\ldots,\mu_{t}^{(i)}-1, \qquad
t=1,2,\ldots,s_{i} \quad {\rm and } \quad i=1,2,\ldots,k. $$
The  coefficients of $P(x)$ can be determined by solving a
linear system of $np^{2}$ unknows and $np^{2}$ equations, because
$\sum_{i=1}^{k} \, \sum_{t=1}^{s_{i}} \, \mu_{t}^{(i)} = np $.

\begin{theorem}
\label{eig:uniopl}
The general polynomial interpolation problem has a unique solution.
\end{theorem}

\noindent
\underline{Proof}\\
Let $P_{1}(x)$ and $P_{2}(x)$ be two solutions of the polynomial interpolation
problem. Then $P_{1}(x) -P_{2}(x)$  is a matrix polynomial of degree $\leq
n-1$ which satisfies
$$ \sum_{l=0}^{q} \frac{1}{l!} \, [P_{1}^{(l)}(x_{i})-P_{2}^{(l)}(x_{i})]
\, v_{t,q-l}^{(i)} \ = \ 0,  $$
where $$ q=0,1,\ldots,\mu_{t}^{(i)}-1, \qquad t=1,2,\ldots,s_{i} \quad
{\rm and} \quad i=1,2,\ldots,k. $$
From this we get that the vectors
$$ v_{t,0}^{(i)}, v_{t,1}^{(i)},\ldots,v_{t,\mu_{t}^{(i)}-1}^{(i)}, \qquad
t=1,2,\ldots,s_{i}, \qquad i=1,2,\ldots,k$$
form a Jordan chain for the matrix polynomial
$P_{1}(x)-P_{2}(x)=B_{n-1}x^{n-1}+ \ldots+B_{1}x+B_{0} $, corresponding to
$x_{i}$. But this implies that $$
B_{n-1}X_{i,t}J_{i,t}^{n-1} +B_{n-2}X_{i,t}J_{i,t}^{n-2} +  \ldots +
B_{1}X_{i,t}J_{i,t} + B_{0}X_{i,t} =0,  $$ where
$ t=1,2,\ldots,s_{i} $ and $ i=1,2,\ldots,k $. The $p
\times \mu_{t}^{(i)}$ dimensional matrix $X_{i,t}$ and the
$ \mu_{t}^{(i)} \times \mu_{t}^{(i)} $ dimensional matrix $J_{i,t}$ are
given by  $$X_{i,t} = \left( \begin{array}{cccc}
 v_{t,0}^{(i)} & v_{t,1}^{(i)} & \ldots & v_{t,\mu_{t}^{(i)}-1}^{(i)} \\
  \end{array}  \right) \quad {\rm and } \quad
 J_{i,t} = \left ( \begin{array}{ccccc}
         x_{i} &   1    &        &        &       \\
               & \ddots & \ddots &        &       \\
               &        & \ddots & \ddots &       \\
               &        &        & x_{i}  &   1    \\
               &        &        &        &  x_{i}
                \end{array} \right ). $$
Hence we get
$$ \left( \begin{array}{cccc}
B_{0} & B_{1} & \ldots & B_{n-1} \\
\end{array} \right) \,
\left( \begin{array}{c}
X_{i,t}  \\   X_{i,t}J_{i,t}  \\  \vdots \\  X_{i,t}J_{i,t}^{n-1}  \\
\end{array} \right) = 0, \qquad  t=1,2,\ldots,s_{i}, \qquad i=1,2,\ldots,k, $$
or
$$ \left( \begin{array}{cccc}
B_{0} & B_{1} & \ldots & B_{n-1} \\
\end{array} \right) \,
\left( \begin{array}{c}
X  \\   XJ  \\  \vdots \\  XJ^{n-1}  \\
\end{array} \right) = 0. $$
Since ${col(XJ^{l})}_{l=0}^{n-1}$ is a non-singular matrix, we have
$$B_{i}=0 \quad  {\rm for} \quad i=1,2,\ldots,n-1. $$  And this implies that
$P_{1}(x)=P_{2}(x)$. \endproof \bigskip

\noindent
\underline{Corollary}\\
If $F(x)$ is a matrix polynomial of degree $ \leq n-1$, then
$P(x)=F(x)$.

The interpolation problem can be formulated in terms of vector polynomials
instead of the sets of $p$-dimensional vectors. Define the $p$-dimensional
vectors $$  \sum_{l=0}^{q} \frac{1}{l!} \, F^{(l)}(x_{i}) \, v_{t,q-l}^{(i)}
= z_{t,q}^{(i)} $$
and the $p$-dimensional vector polynomials
$$ v_{t}^{(i)}(x) = v_{t,0}^{(i)} + v_{t,1}^{(i)} (x-x_{i}) + v_{t,2}^{(i)}
(x-x_{i})^{2} + \ldots + v_{t,\mu_{t}^{(i)}-1}^{(i)} (x-x_{i})^{\mu_{t}^{(i)}-1}
$$
$$ z_{t}^{(i)}(x) = z_{t,0}^{(i)} + z_{t,1}^{(i)} (x-x_{i}) + z_{t,2}^{(i)}
(x-x_{i})^{2} + \ldots + z_{t,\mu_{t}^{(i)}-1}^{(i)} (x-x_{i})^{\mu_{t}^{(i)}-1}
$$ for $ q=0,1,\ldots,\mu_{t}^{(i)}-1 $, $t=1,2,\ldots,s_{i}$ and
$ i=1,2,\ldots,k $.
In this case we are looking for a matrix polynomial  $P(x)$ which satisfies
$$ {\left. \left( \frac{ \partial^{\alpha-1}}{\partial x^{\alpha-1}}\,
P(x)v_{t}^{(i)}(x) \right) \right| }_{x=x_{i}} =
 {\left. \left( \frac{ \partial^{\alpha-1}}{\partial x^{\alpha-1}}\,
z_{t}^{(i)}(x) \right) \right|}_{x=x_{i}} $$
for $\alpha = 1,2,\ldots,\mu_{t}^{(i)}$,
$ t=1,2,\ldots,s_{i}$ and $ i=1,2,\ldots,k$. In addition we know
$v_{t}^{(i)}(x_{i})$ are non-zero and linearly independent vectors. \\
If we set $s_{i}=1$, $i=1,2,\ldots,k$, we have $k$ different points
$x_{1},x_{2},\ldots,x_{k}$, $k$ vector polynomials
$v_{1}^{(1)}(x),v_{1}^{(2)}(x), \ldots, v_{1}^{(k)}(x)$
for which $v_{1}^{(i)}(x_{i}) \not = 0$,  $i=1,2,\ldots,k$ and $k$ vector
polynomials $z_{1}^{(1)}(x),z_{1}^{(2)}(x), \ldots, z_{1}^{(k)}(x)$ and we want
to find a matrix polynomial $P(x)$ which satisfies
$$ {\left. \left( \frac{ \partial^{\alpha-1}}{\partial x^{\alpha-1}}\,
P(x)v_{1}^{(i)}(x) \right) \right|}_{x=x_{i}} =
 {\left.\left( \frac{ \partial^{\alpha-1}}{\partial x^{\alpha-1}}\,
z_{1}^{(i)}(x) \right) \right|}_{x=x_{i}} , \quad \alpha = 1,2,\ldots,m_{i},
\quad i=1,2,\ldots,k. $$  This is the interpolation problem treated in
\cite{baka:1990}.

We will put $\mu_{t}^{(i)}=1$ instead of $s_{i}=1$
in the interpolation problem that we will discuss and which we need for the
construction of the Gaussian quadrature rule.

\subsection{The interpolation problem of Lagrange}
If all $\mu_{j}^{(i)} = 1$, we get the interpolation problem of Lagrange.
This means that we have $k$ different points
$x_{1},x_{2},\ldots,x_{k}$ with  multiplicity resp.\
$m_{1},m_{2},\ldots,m_{k}$ and $\sum_{i=1}^{k} m_{i} = np$.
Every point  $x_{i}$, $i =1,2,\ldots,k$, is associated with a set of non-zero,
linearly independent $p$-dimensional column vectors
$ v_{i,1}, v_{i,2},\ldots,v_{i,m_{i}}$.
Let $X$ and $J$ be defined as follows
$$   X = ( v_{1,1}\ \ldots \ v_{1,m_{1}}\ \ldots \
 v_{k,1}\ \ldots \ v_{k,m_{k}})
\    {\rm and }  \  J =
diag(\underbrace{x_{1},\ldots,x_{1}}_{m_{1}},\underbrace{x_{2},\ldots,x_{2}}_{m_
{2}}, \ldots , \underbrace{x_{k}, \ldots , x_{k}}_{m_{k}} ) , $$
then we require that the matrix $ {col(XJ^{l})}_{l=0}^{n-1} $ is non-singular.

We will construct the interpolation matrix polynomial $P(x)$ of degree
$n-1$ such that  $$ P(x_{i}) \, v_{i,j}  = F(x_{i}) \, v_{i,j},  \qquad
j=1,\ldots,m_{i}, \qquad i=1,\ldots,k. $$

\noindent
\underline{Special case.} \\
If $p=1$, then all the multiplicities have to be equal to $1$.
This means we have $n$ different points $x_{1},x_{2},
\ldots,x_{n} $ and every point $x_{i}$ is associated with a non-zero number
$v_{i}$, $i=1,2,\ldots,k$. The matrix
\begin{eqnarray*}
\left( \begin{array}{c}
       X \\ XJ \\ XJ^{2} \\ \ldots \\ XJ^{n-1}
   \end{array} \right ) & = &
\left( \begin{array}{cccc}
         v_{1}  &      v_{2}    &  \ldots  &    v_{n}  \\
   v_{1} x_{1}  &   v_{2}x_{2}  &  \ldots  & v_{n}x_{n}  \\
   v_{1} x_{1}^{2}  &   v_{2}x_{2}^{2}  &  \ldots  & v_{n}x_{n}^{2}  \\
       \vdots        &    \vdots        &  \vdots  &   \vdots  \\
   v_{1} x_{1}^{n-1} & v_{2}x_{2}^{n-1}  &  \ldots  & v_{n}x_{n}^{n-1}  \\
\end{array} \right )\\
& = &
\left( \begin{array}{cccc}
       1    &    1    &  \ldots  &    1    \\
     x_{1}  &  x_{2}  &  \ldots  &  x_{n}  \\
 x_{1}^{2}  & x_{2}^{2}  &  \ldots  & x_{n}^{2}  \\
       \vdots        &    \vdots        &  \vdots  &   \vdots  \\
 x_{1}^{n-1} &x_{2}^{n-1}  &  \ldots  & x_{n}^{n-1}  \\
\end{array} \right ) \,
\left( \begin{array}{ccccc}
     v_{1}   &    0    &  0       &  \ldots  &    0    \\
     0       &  v_{2}  &  0       &  \ldots  &    0    \\
     0       &    0    &  v_{3}   &  \ldots  &    0    \\
  \vdots     & \vdots  & \vdots   &  \ddots  &   \vdots  \\
     0       &    0    &   0      &  \ldots  &   v_{n} \\
\end{array} \right )
\end{eqnarray*}
is non-singular because of our choice of the points $x_{1},x_{2},\ldots,x_{n}$
and the non zero-numbers $v_{1},v_{2},\ldots,v_{n}$.\\
The interpolation polynomial satisfies
$$ p(x_{i}) v_{i} = f(x_{i}) v_{i}, \qquad  i=1,2,\ldots ,n.$$
But this is equivalent with
 $$ p(x_{i}) = f(x_{i}),  \qquad  i=1,2,\ldots,n,$$
because $v_{i} \not = 0$, for $i=1,2,\ldots,n$.
So we find the polynomial interpolation problem of Lagrange, for real
functions $f : \R \rightarrow \R$.

If $m_{i}=p \cdot d_{i}$, $i=1,2,\ldots,k$ and
${col(X_{i}J_{i}^{l})}_{l=0}^{d_{i}-1}$ are non-singular matrices, then
$(X_{i},J_{i})$ is the Jordan pair of a monic matrix polynomial $L_{i}(x)$ of
degree $d_{i}$. Because $d_{i} \geq 1$, $R_{i}(x) = F(x_{i})$ will be a matrix
polynomial of degree $< d_{i}$ and this for $i=1,2,\ldots,k$. Now we get the
following interpolation problem: given $L_{1}(x), L_{2}(x),\ldots,L_{k}(x)$,
monic matrix polynomials of degree $d_{1},d_{2},\ldots, d_{k}$ respectively and
$R_{1}(x),R_{2}(x),\ldots,R_{k}(x)$ matrix polynomials of degree $< d_{i}$
respectively. Find a matrix polynomial $P(x)$ such that
$$ P(x) = S_{i}(x) L_{i}(x) + R_{i}(x), \qquad i=1,2,\ldots,k $$
for some matrix polynomials $S_{1}(x),S_{2}(x),\ldots,S_{k}(x)$. This
interpolation problem is treated in \cite{goka:1981}, \cite{goka:1982} and
\cite{fu:1983}.

\begin{theorem}
The interpolation matrix polynomial of the Lagrange interpolation problem is
given by
$$ P(x) = \sum_{i=1}^{k} F(x_{i})\, ( 0\  \ldots \ 0 \ v_{i,1}\
\ldots \ v_{i,m_{i}} \ 0\ \ldots \ 0 ) \, (V_{1}+V_{2}x+ \ldots + V_{n}x^{n-1}),
$$ where
$$ (V_{1} \ V_{2}\ \ldots \ V_{n} ) =
{ ( {col(XJ^{l})}_{l=0}^{n-1} )}^{-1},  $$
with $ V_{i}$ a $ np \times p $ matrix.
\end{theorem}

\noindent
\underline{Proof}   \\
The interpolation matrix polynomial $P(x) = A_{n-1}x^{n-1}+\ldots+A_{1}x+A_{0}
$ is determined by a linear system of
$np^{2}$ equations and $np^{2}$ unknows,
$$ P(x_{i}) \, v_{i,j} = F(x_{i}) \, v_{i,j},  \qquad
j=1,\ldots,m_{i}, \qquad i=1,\ldots,k.$$
In block form we get
\begin{eqnarray*}
\lefteqn{
 (A_{0} \ A_{1}\ \ldots \ A_{n-1} )  \cdot
 \left( \begin{array}{c}
       X \\ XJ \\ XJ^{2} \\ \ldots \\ XJ^{n-1}
   \end{array} \right ) = } \\
&& ( F(x_{1})v_{1,1} \  F(x_{1})v_{1,2} \ \ldots \ F(x_{1})v_{1,m_{1}} \
F(x_{2})v_{2,1}\ \ldots \ F(x_{k})v_{k,m_{k}} )
\end{eqnarray*}
(see proof of Theorem \ref{eig:uniopl}). \\
This means that the  coefficients of $P(x)$ satisfy
\begin{eqnarray*}
\lefteqn{ (A_{0} \ A_{1}\ \ldots \ A_{n-1} )\, = } \\
& &  ( F(x_{1})v_{1,1} \  F(x_{1})v_{1,2} \ \ldots \ F(x_{1})v_{1,m_{1}} \
F(x_{2})v_{2,1}\ \ldots \ F(x_{k})v_{k,m_{k}} ) \cdot
 (V_{1} \ V_{2}\ \ldots \ V_{n} )
\end{eqnarray*}
or
$$ A_{j} = ( F(x_{1})v_{1,1} \  F(x_{1})v_{1,2} \ \ldots \ F(x_{1})v_{1,m_{1}} \
F(x_{2})v_{2,1}\ \ldots \ F(x_{k})v_{k,m_{k}} ) \, V_{j+1}, $$
 for $ j=0,1,\ldots,n-1$.
But then we have  $$ P(x) =
 ( F(x_{1})v_{1,1} \ \ldots \ F(x_{1})v_{1,m_{1}} \
F(x_{2})v_{2,1}\ \ldots \ F(x_{k})v_{k,m_{k}} ) \,
(V_{1}+V_{2}x+\ldots+V_{n}x^{n-1}) $$
or
$$ P(x) = \sum_{i=1}^{k} F(x_{i})\, ( 0\  \ldots \ 0 \ v_{i,1}\ \ldots \
v_{i,m_{i}} \ 0\ \ldots \ 0 ) \, (V_{1}+V_{2}x+ \ldots + V_{n}x^{n-1}), $$
with $(m_{1}+\ldots+m_{i-1})$ zeros before $v_{i,1}$ and
$(m_{i+1}+\ldots+m_{k})$ zeros after $v_{i,m_{i}}$.  \endproof

\begin{theorem}
The interpolation matrix polynomial of the Lagrange interpolation problem is
given by
$$ P(x) = \sum_{i=1}^{k} \, F(x_{i})\, W_{i}(x), $$
where
$$ W_{i}(x) = \frac{1}{x-x_{i}}\, \left( \begin{array}{ccc}
 v_{i,1} \ \ldots \ v_{i,m_{i}} \\
 \end{array} \right ) \,
 \left( \begin{array}{c}
  w_{i,1}^{T} \\ \vdots \\  w_{i,m_{i}}^{T} \\
  \end{array}  \right) \, \hat{Q}_{n}(x). $$
The vector $w_{i,j}^{T}$ is the $ (m_{1}+\ldots+m_{i-1}+j)$th row from
$V_{n}$
and  $\hat{Q}_{n}(x)$ is the monic matrix polynomial of degree
$n$ with  Jordan pair $(X,J)$.
\end{theorem}

\noindent
Note that the vectors $w_{i,j}$ are left rootvectors of  the matrix polynomial
$\hat{Q}_{n}(x)$ (Section 1).

\noindent
\underline{Proof} \\
The proof of this theorem consists of three parts :
\begin{itemize}
\item[(a)]  $JV_{1}-V_{n}XJ^{n}V_{1}=0$ and $JV_{i}-V_{n}XJ^{n}V_{i}=V_{i-1}$,
where $i=2,\ldots,n$. \\  \bigskip
To show these relations we consider the similar standard pair
$(X^{\prime},C_{1})$ with
$$ X^{\prime} = ( I \ 0 \ \ldots \ 0 )$$
and $$C_{1} = \left( \begin{array}{ccccc}
   0  &  I  &  0  &  \ldots  &  0  \\
   0  &  0  &  I  &  \ldots  &  0  \\
\vdots &  \vdots &  \vdots  &  \ddots  &  0  \\
   0  &  0  &  0  & \ldots & I \\
  -B_{0} & -B_{1} & -B_{2} & \ldots & -B_{n-1} \\
  \end{array} \right), $$
where $\hat{Q}_{n}(x) = Ix^{n}+B_{n-1}x^{n-1}+\ldots+B_{1}x+B_{0}$, the monic
matrix polynomial with Jordanpair $(X,J)$. This means that $X=X^{\prime}S$ and
$J=S^{-1}C_{1} S$ with $ S = {col(XJ^{l})}_{l=0}^{n-1} $ (see Section 1).
The matrix $ {col(X^{\prime}{C_{1}}^{l})}_{l=0}^{n-1}$ satisfies
$$ \left( \begin{array}{c}
       X ^{\prime} \\ X^{\prime}C_{1} \\ X^{\prime}{C_{1}}^{2} \\ \ldots \\
X^{\prime}{C_{1}}^{n-1} \end{array} \right ) =
   \left( \begin{array}{cccc}
    I  &   0  &  \ldots   & 0 \\
    0  &   I  &  \ldots   & 0 \\
    \vdots & \vdots & \ddots & \vdots \\
    0  &   0  &  \ldots   & I \\
\end{array} \right),$$
because
\begin{eqnarray*}
\lefteqn{ X^{\prime}{C_{1}}^{k} } \\
& = & ( 0 \ \ldots \ I \ \ldots \ 0 )  \quad {\rm with\
} I \ {\rm on\ the \ (k+1)th\ place\ and \ } k=0,1,\ldots,n-1,\\
& = & ( -B_{0}\  -B_{1}\ \ldots \ -B_{n-1})\quad {\rm for  } \quad k=n.
\end{eqnarray*}
So the matrices $V_{i}^{\prime}$ are given by
$$ V_{i}^{\prime} = \left( \begin{array}{c}
0  \\ \vdots \\ I \\ \vdots  \\ 0 \\
\end{array} \right) \quad {\rm with\ } I \ {\rm on\ the \ }i{\rm th\ place}.$$
But then we have
$$ C_{1}\,
\left( \begin{array}{c}
I  \\ \vdots \\ 0 \\ \vdots  \\ 0 \\
\end{array} \right) -
\left( \begin{array}{c}
0  \\ \vdots \\ 0 \\ \vdots  \\ I \\
\end{array} \right) \, ( -B_{0}\  -B_{1}\ \ldots \ -B_{n-1}) \,
\left( \begin{array}{c}
I  \\ \vdots \\ 0 \\ \vdots  \\ 0 \\
\end{array} \right)   =
\left( \begin{array}{c}
0  \\ \vdots \\ 0 \\ \vdots  \\ -B_{0} \\
\end{array} \right)  -
\left( \begin{array}{c}
0  \\ \vdots \\ 0 \\ \vdots  \\ -B_{0} \\
\end{array} \right)         =
\left( \begin{array}{c}
0  \\ \vdots \\ 0 \\ \vdots  \\ 0 \\
\end{array} \right) $$
or
$$ C_{1} V_{1}^{\prime}-V_{n}^{\prime} X^{\prime} C_{1}^{n} V_{1}^{\prime} = 0.
$$
Consequently,
$$ S^{-1}C_{1}SS^{-1} V_{1}^{\prime} -S^{-1} V_{n}^{\prime} X^{\prime}SS^{-1}
C_{1}^{n}SS^{-1} V_{1}^{\prime} = 0 $$
or         $$JV_{1}-V_{n}XJ^{n}V_{1}=0.$$
For $i=2,3,\ldots,n$ we have
\begin{eqnarray*}
\lefteqn{
 C_{1}\,
\left( \begin{array}{c}
0  \\ \vdots \\ 0 \\ I \\ \vdots  \\ 0 \\
\end{array} \right)
-
\left( \begin{array}{c}
0  \\ \vdots \\ 0  \\ 0 \\ \vdots  \\ I \\
\end{array} \right) \, ( -B_{0}\  -B_{1}\ \ldots \ -B_{n-1}) \,
\left( \begin{array}{c}
0  \\ \vdots \\ 0 \\ I \\ \vdots  \\ 0 \\
\end{array} \right) } \\
& = &
\left( \begin{array}{c}
0  \\ \vdots \\ I \\ 0 \\ \vdots  \\ -B_{i-1}\\
\end{array} \right) \, - \,
\left( \begin{array}{c}
0  \\ \vdots \\ 0 \\ 0 \\ \vdots  \\ -B_{i-1} \\
\end{array} \right)        \, = \,
\left( \begin{array}{c}
0  \\ \vdots \\ I \\ 0 \\ \vdots  \\ 0 \\
\end{array} \right).
\end{eqnarray*}
In other words
$$ C_{1} V_{i}^{\prime}-V_{n}^{\prime} X^{\prime} C_{1}^{n} V_{i}^{\prime} =
V_{i-1}^{\prime} $$
or
$$ S^{-1}C_{1}SS^{-1} V_{i}^{\prime} -S^{-1} V_{n}^{\prime} X^{\prime}SS^{-1}
C_{1}^{n}SS^{-1} V_{i}^{\prime} = S^{-1} V_{i-1}^{\prime}, $$
so that
         $$JV_{i}-V_{n}XJ^{n}V_{i}=V_{i-1}. $$

\item[(b)] $V_{n}\hat{Q}_{n}(x) = (xI-J)\, (V_{1}+V_{2}x+\ldots+V_{n}x^{n-1})$
\\ \bigskip
The monic matrix polynomial  $\hat{Q}_{n}(x)$ with Jordan pair $(X,J)$
is given by
 $$\hat{Q}_{n}(x) = x^{n}I - XJ^{n} (V_{1}+V_{2}x+\ldots+V_{n}x^{n-1}). $$
(see Section 1).
But this implies we immediately get the required expression from the relations
proved in (a).

\item[(c)] $ W_{i}(x) = \frac{1}{x-x_{i}}\, \left( \begin{array}{ccc}
 v_{i,1} \ \ldots \ v_{i,m_{i}} \\
 \end{array} \right ) \,
 \left( \begin{array}{c}
  w_{i,1}^{T} \\ \vdots \\  w_{i,m_{i}}^{T} \\
  \end{array}  \right) \, \hat{Q}_{n}(x). $ \\ \bigskip
From the property proved in (b) we get
$$ {(xI-J)}^{-1} V_{n} \hat{Q}_{n}(x) = (V_{1}+V_{2}x+\ldots+V_{n}x^{n-1}).$$
If we call $I_{m_{i}}$ the $np \times np $ diagonal matrix with  $1$ on
the $(m_{1}+\ldots+m_{i-1}+j)$th
row, where $j=1,2,\ldots,m_{i}$ and $0$ on all the other places
and we multiply both sides of the equation on the left with this matrix,
we get
$$ \frac{1}{x-x_{i}}\,
 \left( \begin{array}{c}
0\\ \vdots \\  w_{i,1}^{T} \\ \vdots \\  w_{i,m_{i}}^{T} \\ \vdots \\ 0 \\
  \end{array}  \right) \, \hat{Q}_{n}(x) =  I_{m_{i}} \,
  (V_{1}+V_{2}x+\ldots+V_{n}x^{n-1}).$$
Multiplication on the left of this equation with
the matrix $X$, gives
\begin{eqnarray*}
\lefteqn{ \frac{1}{x-x_{i}}\,
 (  v_{i,1}\ v_{i,2} \ \ldots \ v_{i,m_{i}}) \,
 \left( \begin{array}{c}
  w_{i,1}^{T} \\ w_{i,2}^{T} \\  \vdots \\  w_{i,m_{i}}^{T} \\
  \end{array}  \right) \, \hat{Q}_{n}(x) = }\\
& &   ( 0\ 0\ \ldots \ 0 \ v_{i,1}\
\ldots \ v_{i,m_{i}} \ 0\ \ldots \ 0)\, (V_{1}+V_{2}x+ \ldots + V_{n}x^{n-1}),
\end{eqnarray*}
such that the previous theorem gives the required expression
$$ P(x) = \sum_{i=1}^{k} \, F(x_{i})\,
 \frac{1}{x-x_{i}}\, \left( \begin{array}{ccc}
 v_{i,1} \ \ldots \ v_{i,m_{i}} \\
 \end{array} \right ) \,
 \left( \begin{array}{c}
  w_{i,1}^{T} \\ \vdots \\  w_{i,m_{i}}^{T} \\
  \end{array}  \right) \, \hat{Q}_{n}(x). \ \ \rule{0.5em}{0.5em}$$
\end{itemize}
\bigskip

\noindent
We get another representation of the interpolation matrix polynomial if we
choose the points $x_{i}$ and the corresponding vectors in a special way.
This representation shall be very useful for the construction of quadrature
formulas.

\begin{theorem}
\label{eig:intlagr}
Let $(X,J)$ be a Jordan pair of the orthonormal matrix polynomial
$P_{n}(x)$, then
the interpolation matrix polynomial of the Lagrange interpolation problem is
given by
$$ P(x) = \sum_{i=1}^{k} \, F(x_{i})\, W_{i}(x),$$
with $k$ the number of different zeros $x_{i}$,
$m_{i}$
the multiplicity of $x_{i}$, $v_{i,j}$ the vectors associated with
$x_{i}$ and $$ W_{i}(x) = \left( \begin{array}{ccc}
 v_{i,1} \ \ldots \ v_{i,m_{i}} \\
 \end{array} \right ) \, {K_{i}}^{-1} \,
 \left( \begin{array}{c}
  v_{i,1}^{T} \\ \vdots \\  v_{i,m_{i}}^{T} \\
  \end{array}  \right) \, \frac{P_{n+1}(x_{i})^{T}D_{n+1}P_{n}(x)}{x-x_{i}}. $$
The $m_{i} \times m_{i}$ dimensional matrix $K_{i}$ satisfies
$$  K_{i}  =   -\,
 \left( \begin{array}{c}
  v_{i,1}^{T} \\ \vdots \\  v_{i,m_{i}}^{T} \\
  \end{array}  \right) \, K_{n-1}(x_{i},x_{i}) \,
 \left( \begin{array}{ccc}
 v_{i,1} \ \ldots \ v_{i,m_{i}} \\
 \end{array} \right ), $$
with $$ K_{n-1}(x,y) = \sum_{j=0}^{n-1} P_{j}(y)^{T} P_{j}(x). $$
\end{theorem}

\noindent
Note that the Jordan chains of the orthonormal matrix polynomials $P_{n}(x)$
have length $1$. So, if  $x_{1},x_{2}, \ldots,x_{k}$ are the zeros of the
orthonormal matrix polynomials $P_{n}(x)$ with multiplicity resp.\
$m_{1},m_{2},\ldots,m_{k}$, then every $x_{i}$ is associated with
non-zero and linearly independent vectors $v_{i,1},\ldots,v_{i,m_{i}}$.
Furthermore $\sum_{i=1}^{k} m_{i}=np$ and the matrix
${col(XJ^{l})}_{l=0}^{n-1}$ is non-singular.
In other words, these are exactly the requirements which the points and the
vectors for the Lagrange interpolation problem, have to fulfil.

\noindent
\underline{Proof} \\
From the identity of Christoffel-Darboux and from the fact that $v_{i,j}$
is a rootvector of $P_{n}(x)$ associated with the zero $x_{i}$, we get
\begin{equation}
\label{eq:rvch}
 v_{i,j}^{T} \,\sum_{l=0}^{n-1} P_{l}(x_{i})^{T} P_{l}(x) =
 - \, v_{i,j}^{T} \frac{P_{n+1}(x_{i})^{T}D_{n+1}P_{n}(x)}{x-x_{i}}.
\end{equation}
This implies that the left side of the equation, as well the right side is a
vector polynomial of degree $n-1$.
In other words the $m_{i} \times p $ dimensional matrix
$$ \left( \begin{array}{c}
  v_{i,1}^{T} \\ \vdots \\  v_{i,m_{i}}^{T} \\
  \end{array}  \right) \,
P_{n+1}(x_{i})^{T}D_{n+1}P_{n}(x) $$ is divisible by $x-x_{i}$.\\
Set
$$ W_{i}(x) = A \,
\left( \begin{array}{c}
  v_{i,1}^{T} \\ \vdots \\  v_{i,m_{}i}^{T} \\
  \end{array}  \right) \,
\frac{P_{n+1}(x_{i})^{T}D_{n+1}P_{n}(x)}{x-x_{i}},
$$ then this matrix polynomial of degree $n-1$ satisfies
$$ W_{i}(x_{l})\, v_{l,j} = 0, \qquad  j=1,2,\ldots,m_{l}, \quad
l=1,\ldots ,i-1,i+1, \ldots k, $$
because there is no singularity for $x_{l} \not =  x_{i}$
and $P_{n}(x_{l})v_{l,j}=0$.
Now we determine the $p \times m_{i}$ dimensional matrix $A$ such that
$$ W_{i}(x_{i})\, v_{i,j} = v_{i,j}, \qquad j=1,2,\ldots,m_{i}. $$
From (\ref{eq:rvch}) we get the following equation :
\begin{eqnarray*}
\lefteqn{ v_{i,j}^{T} \,\sum_{l=0}^{n-1} P_{l}(x_{i})^{T} P_{l}(x)\, v_{i,s} =
}\\
& &  - \, \left( v_{i,j}^{T} \frac{P_{n+1}(x_{i})^{T}D_{n+1}P_{n}(x)}{x-x_{i}}
v_{i,s}
- v_{i,j}^{T} \frac{P_{n+1}(x_{i})^{T}D_{n+1}P_{n}(x_{i})}{x-x_{i}} v_{i,s}
\right). \end{eqnarray*}
Let $x$ approach $x_{i}$, then this becomes
\begin{eqnarray*}
 v_{i,j}^{T} \,\sum_{l=0}^{n-1} P_{l}(x_{i})^{T} P_{l}(x_{i})\, v_{i,s} & = &
 - \, v_{i,j}^{T} P_{n+1}(x_{i})^{T} D_{n+1} P_{n}^{\prime}(x_{i}) v_{i,s} \\
 & = & - {(K_{i})}_{j,s}.
 \end{eqnarray*}
So, the matrix $A$ satisfies
$$ A \,
\left( \begin{array}{c}
  v_{i,1}^{T} \\ \vdots \\  v_{i,m_{i}}^{T} \\
  \end{array}  \right) \, P_{n+1}(x_{i})^{T}D_{n+1}P_{n}^{\prime}(x_{i}) \,
v_{i,s} = v_{i,s}, \qquad s=1,2,\ldots,m_{i}, $$
or
$$ A \,
\left( \begin{array}{c}
  {(K_{i})}_{1,s} \\ {(K_{i})}_{2,s} \\ \vdots \\ {(K_{i})}_{m_{i},s} \\
  \end{array}  \right) = v_{i,s}, \qquad  s=1,2,\ldots,m_{i}. $$
Consequently,  $A$ satisfies
$$ A \,
\left( \begin{array}{cccc}
  {(K_{i})}_{1,1} & {(K_{i})}_{1,2} & \ldots & {(K_{i})}_{1,m_{i}} \\
  {(K_{i})}_{2,1} & {(K_{i})}_{2,2} & \ldots & {(K_{i})}_{2,m_{i}} \\
    \vdots        &     \vdots      &  \vdots &   \vdots   \\
{(K_{i})}_{m_{i},1} & {(K_{i})}_{m_{i},2} & \ldots & {(K_{i})}_{m_{i},m_{i}} \\
\end{array}  \right)\, = \,
\left( \begin{array}{cccc}
 v_{i,1} & v_{i,2} & \ldots & v_{i,m_{i}} \\
 \end{array} \right ) $$
or
$$ A =
\left( \begin{array}{cccc}
 v_{i,1} & v_{i,2} & \ldots & v_{i,m_{i}} \\
 \end{array} \right ) \, {K_{i}}^{-1}. $$
Bringing everything together, we get
$$ W_{i}(x) = \left( \begin{array}{ccc}
 v_{i,1} \ \ldots \ v_{i,m_{i}} \\
 \end{array} \right ) \, {K_{i}}^{-1} \,
 \left( \begin{array}{c}
  v_{i,1}^{T} \\ \vdots \\  v_{i,m_{i}}^{T} \\
  \end{array}  \right) \, \frac{P_{n+1}(x_{i})^{T}D_{n+1}P_{n}(x)}{x-x_{i}}, $$
where $$
 K_{i}  =   -\,
 \left( \begin{array}{c}
  v_{i,1}^{T} \\ \vdots \\  v_{i,m_{i}}^{T} \\
  \end{array}  \right) \, K_{n-1}(x_{i},x_{i}) \,
 \left( \begin{array}{ccc}
 v_{i,1} \ \ldots \ v_{i,m_{i}} \\
 \end{array} \right ) .  $$
Finally we have to show that the matrix $K_{i}$ is non-singular. Therefore we
show that $-K_{i}$ is a symmetric and positive definite matrix.
From the fact that $K_{n-1}(x_{i},x_{i})$ is a symmetric and
positive definite matrix, we get
$$ - K_{i}^{T}  =
 \left( \begin{array}{c}
  v_{i,1}^{T} \\ \vdots \\  v_{i,m_{i}}^{T} \\
  \end{array}  \right) \, K_{n-1}(x_{i},x_{i}) \,
 \left( \begin{array}{ccc}
 v_{i,1} \ \ldots \ v_{i,m_{i}} \\
 \end{array} \right )  = - K_{i}.  $$
Let $w$ be a  $m_{i}$ dimensional vector, then we have
\begin{eqnarray*}
\lefteqn{
 \left( \begin{array}{ccc} w_{1} & \ldots &  w_{m_{i}} \\ \end{array}  \right)
  \, (-K_{i}) \,
 \left( \begin{array}{c} w_{1} \\ \vdots \\  w_{m_{i}} \\ \end{array}  \right)
  } \\
& = &
 \left( \begin{array}{ccc} w_{1} & \ldots &  w_{m_{i}} \\ \end{array}  \right)
\, \left( \begin{array}{c} v_{i,1}^{T} \\ \vdots \\  v_{i,m_{i}}^{T} \\
 \end{array}  \right)
  \, K_{n-1}(x_{i},x_{i}) \,
\left( \begin{array}{ccc} v_{i,1} \ \ldots \ v_{i,m_{i}} \\ \end{array} \right)
\left( \begin{array}{c} w_{1} \\ \vdots \\  w_{m_{i}} \\ \end{array}  \right) \\
& = &
( w_{1}v_{i,1}^{T}\, +\, \ldots\, +\,  w_{m_{i}}v_{i,m_{i}}^{T})
 \, K_{n-1}(x_{i},x_{i}) \,
( w_{1}v_{i,1}\, +\, \ldots\, +\,  w_{m_{i}}v_{i,m_{i}}).
\end{eqnarray*}
This expression is always $\geq 0 $ and only equal to $0$ if the vector
$$  w_{1}v_{i,1}\, +\, \ldots\, +\,  w_{m_{i}}v_{i,m_{i}} $$ is equal to $0$.
But the vectors $v_{i,1},\ldots,v_{i,m_{i}}$ are linearly independent, so that
$ w^{T}(-K_{i}) w = 0 $ if and only if $w=0$.
So $-K_{i}$ is a symmetric and positive definite matrix and taking the inverse
of it will not be a problem. \endproof
\bigskip

\noindent
\underline{Special case.} \\
If $p=1$ we consider $n$ different points $x_{1},\ldots,x_{n}$
and associate a non zero number $v_{i}$ with every point $x_{i}$.
The interpolation polynomial for a function $f(x)$ is given by
 $$ p(x) = \sum_{i=1}^{n} f(x_{i})\,
w_{i}(x), $$ with
\begin{eqnarray*}
 w_{i}(x) &  = & - v_{i} \frac{-1}{v_{i}K_{n-1}(x_{i},x_{i})v_{i}} v_{i}
\frac{p_{n+1}(x_{i}) d_{n+1} p_{n}(x)}{x-x_{i}} \\
& = & \frac{-p_{n+1}(x_{i}) k_{n,n}}{k_{n+1,n+1} K_{n-1}(x_{i},x_{i})}
\frac{ p_{n}(x)}{x-x_{i}} \\
& = & \frac{ p_{n}(x)}{(x-x_{i}) p_{n}^{\prime}(x_{i})}. \\
\end{eqnarray*}

\noindent
\underline{Example.} \\
Consider the orthonormal matrix polynomials on the
interval $[-1,1]$ with respect to the weight matrix function
 $$ W(x) = \left( \begin{array}{cc}
\frac{1}{\pi} {(1-x^{2})}^{-1/2}  &     0   \\
 0   & \frac{2}{\pi} {(1-x^{2})}^{1/2}    \\
\end{array} \right ). $$
The recursion coefficients are given by :
$$ E_{n} = 0 \quad {\rm for} \quad n \geq 0, $$
$$ D_{1} =
 \left( \begin{array}{cc}
\frac{1}{\sqrt{2}}  &     0   \\
0   & \frac{1}{2}    \\
\end{array} \right ) \quad {\rm and} \quad
  D_{n} =
\left( \begin{array}{cc}
\frac{1}{2}  &     0   \\
    0        & \frac{1}{2}    \\
\end{array} \right )  \quad {\rm for} \quad  n \geq 1. $$
In this case the orthonormal matrix polynomials are
$$ P_{0}(x) =
 \left( \begin{array}{cc}
 1  & 0   \\
 0  & 1   \\
\end{array} \right ), \quad
P_{1}(x) =
 \left( \begin{array}{cc}
 \sqrt{2}x  & 0   \\
 0  & 2x   \\
\end{array} \right ),$$
$$ P_{2}(x) =
 \left( \begin{array}{cc}
 2 \sqrt{2} x^{2} - \sqrt{2}  & 0   \\
 0  & 4 x^{2}-1   \\
\end{array} \right ), \quad
P_{3}(x) =
 \left( \begin{array}{cc}
 4 \sqrt{2} x^{3} - 3 \sqrt{2}x  & 0   \\
 0  & 8 x^{3}-4x  \\
\end{array} \right ).$$
In general, we have
$$ P_{n}(x) =
 \left( \begin{array}{cc}
  \sqrt{2}\, T_{n}(x)  & 0   \\
 0  & U_{n}(x)   \\
\end{array} \right ). $$
Consider $P_{2}(x)$, then the zeros and the corresponding rootvectors are given
by
$$ x_{1}=\frac{1}{\sqrt{2}}, \quad \quad x_{2}=\frac{-1}{\sqrt{2}}, \quad \quad
x_{3}=\frac{1}{2}, \quad \quad x_{4}=\frac{-1}{2}, $$
$$ v_{1}= \left( \begin{array}{c}  1  \\  0 \end{array} \right), \quad
   v_{2}= \left( \begin{array}{c}  4  \\  0 \end{array} \right), \quad
   v_{3}= \left( \begin{array}{c}  0  \\  3 \end{array} \right), \quad
   v_{4}= \left( \begin{array}{c}  0  \\ -2 \end{array} \right).$$
This leads to
$$ K_{1} = - \left( \begin{array}{cc} 1 & 0   \end{array} \right) \,
\left( \begin{array}{cc}  1+1 & 0 \\ 0 & 1+2 \\ \end{array} \right) \,
 \left( \begin{array}{c}  1  \\  0 \end{array} \right) = -2 $$
and
\begin{eqnarray*}
 W_{1}(x) &  = &
 \left( \begin{array}{c}  1  \\  0 \end{array} \right)  \frac{-1}{2}
 \left( \begin{array}{cc} 1 & 0   \end{array} \right)
\frac{1}{x-\frac{1}{\sqrt{2}}}
\left( \begin{array}{cc}  -1 & 0 \\ 0 & 0 \\ \end{array} \right)
\left( \begin{array}{cc}  \frac{1}{2} & 0 \\ 0 & \frac{1}{2} \\
\end{array} \right) \left( \begin{array}{cc}
 2 \sqrt{2} x^{2} - \sqrt{2}  & 0   \\
 0  & 4 x^{2}-1   \\
\end{array} \right ) \\
& = &
 \left( \begin{array}{cc}
 \frac{1}{\sqrt{2}} (x + \frac{1}{\sqrt{2}})  & 0   \\
 0  & 0   \\
\end{array} \right ).
\end{eqnarray*}
In a similar manner we find
$$ K_{2} = -32, \qquad  W_{2}(x)=
 \left( \begin{array}{cc}
 \frac{-1}{\sqrt{2}} (x - \frac{1}{\sqrt{2}})  & 0   \\
 0  & 0   \\
\end{array} \right ),     $$
$$ K_{3} = -18, \qquad W_{3}(x) =
 \left( \begin{array}{cc}
 0  & 0   \\
0 & (x + \frac{1}{2})   \\
\end{array} \right ),     $$
$$ K_{4} = -8, \qquad W_{4}(x) =
 \left( \begin{array}{cc}
 0  & 0   \\
0 & -( x - \frac{1}{2})   \\
\end{array} \right ).     $$
Let $F(x)$ be given by
$$F(x) =
 \left( \begin{array}{cc}
 2x+5  & 6x   \\
   7   & 4x-3  \\
\end{array} \right ), $$
then the interpolating matrix polynomial is given by
\begin{eqnarray*}
P(x) & = &
 \left( \begin{array}{cc}
 5+\sqrt{2}  & 3\sqrt{2}   \\
   7   & 2 \sqrt{2}-3  \\
\end{array} \right ) \, W_{1}(x) \, + \,
 \left( \begin{array}{cc}
 5-\sqrt{2}  & -3\sqrt{2}   \\
   7   & -2 \sqrt{2}-3  \\
\end{array} \right ) \, W_{2}(x) \\
& + &
 \left( \begin{array}{cc}
 6  & 3   \\
  7 & -1  \\
\end{array} \right ) \, W_{3}(x) \, + \,
 \left( \begin{array}{cc}
 4  & -3   \\
  7 & -5  \\
\end{array} \right ) \, W_{4}(x) \\
& = & \left( \begin{array}{cc}
 2x+5  & 6x   \\
   7   & 4x-3  \\
\end{array} \right ) = F(x).
\end{eqnarray*}
Let $F(x)$ be given by
$$F(x) =
 \left( \begin{array}{cc}
 x^{2}+1  & 6x   \\
 7x+1    & 5x^{2}-1  \\
\end{array} \right ), $$
then the interpolating matrix polynomial is given by
\begin{eqnarray*}
P(x) & = &
 \left( \begin{array}{cc}
 \frac{3}{2}  & 3\sqrt{2}   \\
 \frac{7}{\sqrt{2}}+1  & \frac{3}{2}   \\
\end{array} \right ) \, W_{1}(x) \, + \,
 \left( \begin{array}{cc}
 \frac{3}{2}  & -3\sqrt{2}   \\
\frac{-7}{\sqrt{2}}+1  & \frac{3}{2}   \\
\end{array} \right ) \, W_{2}(x) \\
& + & \left( \begin{array}{cc}
\frac{5}{4}  & 3   \\
\frac{9}{2} & \frac{1}{4}  \\
\end{array} \right ) \, W_{3}(x) \, + \,
 \left( \begin{array}{cc}
\frac{5}{4}  & -3   \\
\frac{-5}{2} & \frac{1}{4}  \\
\end{array} \right ) \, W_{4}(x) \\
& = & \left( \begin{array}{cc}
\frac{3}{2}  & 6x   \\
   7x+1   & \frac{1}{4}  \\
\end{array} \right ).
\end{eqnarray*}

\section{Gaussian quadrature}

\noindent
Let $W(x)$ be a  matrix weight function defined on the interval
$[a,b]$. Then we are going to approximate the integral of
matrix functions by means of a sum of the form
 $$ \int_{a}^{b} F(x)\, W(x)\, G(x)^{T}\, dx  \simeq \sum_{i=1}^{k} F(x_{i}) \,
\Lambda_{i} \, G(x_{i})^{T}, $$
where $ \Lambda_{i} \in \Rpp$.

\subsection{Quadrature formulas}
The numerical integration problem exists in finding points
$x_{i}$ and matrices $\Lambda_{i}$, where $i=1,\ldots,k$, so that the formula
will have the highest possible degree of precision.
This degree of precision of a quadrature formula is equal to $m (\in \N)$ if
\begin{eqnarray*}
\int_{a}^{b} x^{l} \, W(x)\, dx - \sum_{i=1}^{k} x_{i}^{l}\,
\Lambda_{i} & = & 0 \quad \quad {\rm for} \quad l=0,1,\ldots,m, \\
 & \not  = 0 & \quad \quad {\rm for} \quad l=m+1.
\end{eqnarray*}
Note that this agrees with the accuracy of the quadrature formula
$$ \int_{a}^{b} F(x)\, W(x)\, G(x)^{T}\, dx \simeq \sum_{i=1}^{k} F(x_{i}) \,
\Lambda_{i}\, G(x_{i})^{T},$$
for matrix polynomials which satisfy  $\deg{F(x)}+\deg{G(x)} \leq m$. \\
It will be convient to choose
$$ \Lambda_{i} =
\left( \begin{array}{ccc}
 v_{i,1} &  \ldots & v_{i,m_{i}} \\
 \end{array} \right ) \, {A_{i}} \,
\left( \begin{array}{c}
  v_{i,1}^{T} \\ \vdots \\  v_{i,m_{i}}^{T} \\
  \end{array}  \right), $$
with $v_{i,1}, \ldots, v_{i,m_{i}}$ linearly independent, non-zero vectors and
$\sum_{i=1}^{k} m_{i} = np $.
The numerical integration problem exists than in determining
$x_{i}$, vectors $v_{i,j}$ and matrices $A_{i}$, where
$j=1,\ldots,m_{i}$ and $i=1,\ldots,k$, so that the formula will have the
highest possible degree of precision.
We will show now that $2n-1$ is the highest degree of precision
for a quadrature formula of the postulated form.

\begin{theorem}
A quadrature formula of the form
$$ \int_{a}^{b} F(x) \, W(x) \, G(x)^{T}\, dx  \simeq \sum_{i=1}^{k} F(x_{i})\,
\Lambda_{i}\, G(x_{i})^{T} $$
cannot be exact for all matrix polynomials $F(x)$ and $G(x)$ satisfying
$$ \deg{F(x)}+\deg{G(x)} \leq 2n.$$
\end{theorem}

\noindent
\underline{Proof}\\
Let $F(x) = G(x) = \hat{Q}_{n}(x) $, where $\hat{Q}_{n}(x)$ is the monic
matrix polynomial of degree $n$ which satisfies
$\hat{Q}_{n}(x_{i})v_{i,j} = 0 $ for $j=1,2,\ldots,m_{i}$ and $i=1,2,\ldots,k$.
In this case we have
$$ \int_{a}^{b} F(x) \, W(x) \, G(x)^{T}\, dx  =
 \int_{a}^{b} \hat{Q}_{n}(x) \, W(x) \, {\hat{Q}_{n}(x)}^{T}\, dx  \, > 0. $$
On the other hand we have
$$ \sum_{i=1}^{k} F(x_{i})\, \Lambda_{i}\, G(x_{i})^{T}  =
 \sum_{i=1}^{k} \hat{Q}_{n}(x_{i})\,
\left( \begin{array}{ccc}
 v_{i,1} &  \ldots & v_{i,m_{i}} \\
 \end{array} \right ) \, {A_{i}} \,
\left( \begin{array}{c}
  v_{i,1}^{T} \\ \vdots \\  v_{i,m_{i}}^{T} \\
  \end{array}  \right)   \, {\hat{Q}_{n}(x_{i})}^{T} \, = 0, $$
so that the quadrature formula cannot be exact in this case.
\endproof\bigskip

\noindent
In the following properties we will show that it is possible to construct a
quadrature formula with degree of precision $2n-1$.

\begin{theorem}
Let $(X,J)$ be a Jordan pair of the orthonormal matrix polynomial
$P_{n}(x)$ on the interval $[a,b]$ and with respect to the weight matrix
function $W(x)$. Then we have
\begin{equation}
\label{eq:gakwa}
 \int_{a}^{b} F(x) \, W(x) \, G(x)^{T}\, dx  \simeq \sum_{i=1}^{k} F(x_{i})\,
\Lambda_{i}\, G(x_{i})^{T},
\end{equation}
where $k$ is the number of different zeros $x_{i}$, $m_{i}$ is the
multiplicity of $x_{i}$, $v_{i,j}$ are the vectors associated with $x_{i}$,
$$ \Lambda_{i} =
\left( \begin{array}{ccc}
 v_{i,1} &  \ldots & v_{i,m_{i}} \\
 \end{array} \right ) \, {L_{i}}^{-1} \,
\left( \begin{array}{c}
  v_{i,1}^{T} \\ \vdots \\  v_{i,m_{i}}^{T} \\
  \end{array}  \right)  $$
and
 $$ L_{i} =
\left( \begin{array}{c}
  v_{i,1}^{T} \\ \vdots \\  v_{i,m_{i}}^{T} \\
  \end{array}  \right)  \, K_{n-1}(x_{i},x_{i}) \,
\left( \begin{array}{ccc}
 v_{i,1} & \ldots & v_{i,m_{i}} \\
 \end{array} \right ). $$
This quadrature formula is exact for matrix polynomials
$F(x)$ and $G(x)$ which satisfy  $$ \deg{F(x)}+\deg{G(x)} \leq 2n-1.$$
\end{theorem}

\noindent
\underline{Proof} \\
If we replace $F(x)$ and $G(x)$ by their interpolating matrix polynomial from
the Lagrange interpolation problem (see Theorem \ref{eig:intlagr}), then
we get  $$ \int_{a}^{b} F(x) \, W(x) \, G(x)^{T}\, dx  \simeq
\sum_{i=1}^{k} \sum_{j=1}^{k} F(x_{i})\, \Lambda_{i,j} \, \, G(x_{j})^{T}, $$
with
\begin{eqnarray*}
\Lambda_{i,j} & =  &
\left( \begin{array}{ccc}
 v_{i,1} &  \ldots & v_{i,m_{i}} \\
 \end{array} \right ) \, {K_{i}}^{-1} \,
\left( \begin{array}{c}
  v_{i,1}^{T} \\ \vdots \\  v_{i,m_{i}}^{T} \\
  \end{array}  \right) \, \\
& & \times {P_{n+1}(x_{i})}^{T}\, D_{n+1}\, \int_{a}^{b}
\frac{P_{n}(x)}{x-x_{i}}\, W(x)\, \frac{P_{n}(x)^{T}}{x-x_{j}}\, dx \,
 D_{n+1}\, P_{n+1}(x_{j})  \\
& &  \times \left( \begin{array}{ccc}
 v_{j,1} &  \ldots & v_{j,m_{i}} \\
 \end{array} \right ) \, {K_{j}}^{-1} \,
\left( \begin{array}{c}
  v_{j,1}^{T} \\ \vdots \\  v_{j,m_{j}}^{T} \\
  \end{array}  \right).
\end{eqnarray*}
From the identity of Christoffel-Darboux we get
$$ \left( \begin{array}{c}
  v_{i,1}^{T} \\ \vdots \\  v_{i,m_{i}}^{T} \\
  \end{array}  \right) \, K_{n-1}(x,x_{i}) =
- \left( \begin{array}{c}
  v_{i,1}^{T} \\ \vdots \\  v_{i,m_{i}}^{T} \\
  \end{array}  \right) \, \frac{P_{n+1}(x_{i})^{T}D_{n+1}P_{n}(x)}{x-x_{i}},
$$ so that
\begin{eqnarray*}
\Lambda_{i,j} & = &
\left( \begin{array}{ccc}
 v_{i,1} &  \ldots & v_{i,m_{i}} \\
 \end{array} \right ) \, {K_{i}}^{-1} \,
\left( \begin{array}{c}
  v_{i,1}^{T} \\ \vdots \\  v_{i,m_{i}}^{T} \\
  \end{array}  \right) \,
 \int_{a}^{b} K_{n-1}(x,x_{i})\, W(x)\, {K_{n-1}(x,x_{j})}^{T}\, dx \,  \\
&   &   \times \, \left( \begin{array}{ccc}
 v_{j,1} &  \ldots & v_{j,m_{i}} \\
 \end{array} \right ) \, {K_{j}}^{-1} \,
\left( \begin{array}{c}
  v_{j,1}^{T} \\ \vdots \\  v_{j,m_{j}}^{T} \\
  \end{array}  \right) \\
& = &
\left( \begin{array}{ccc}
 v_{i,1} &  \ldots & v_{i,m_{i}} \\
 \end{array} \right ) \, {K_{i}}^{-1}
\left( \begin{array}{c}
  v_{i,1}^{T} \\ \vdots \\  v_{i,m_{i}}^{T} \\
  \end{array}  \right)
 K_{n-1}(x_{j},x_{i})
\left( \begin{array}{ccc}
 v_{j,1} &  \ldots & v_{j,m_{i}} \\
 \end{array} \right ) \, {K_{j}}^{-1}
\left( \begin{array}{c}
  v_{j,1}^{T} \\ \vdots \\  v_{j,m_{j}}^{T} \\
  \end{array}  \right).
\end{eqnarray*}
But for  $i \not = j$ we have
\begin{eqnarray*}
v_{i,s}^{T} K_{n-1}(x_{j},x_{i}) v_{j,t}  & = &
v_{i,s}^{T} K_{n}(x_{j},x_{i}) v_{j,t} \\
& = & v_{i,s}^{T}
 \frac{P_{n}(x_{i})^{T}D_{n+1}P_{n+1}(x_{j})-
P_{n+1}(x_{i})^{T}D_{n+1}P_{n}(x_{j})}{x_{j}-x_{i}} v_{j,t} \\
& = & 0,
\end{eqnarray*}
so that  $\Lambda_{i,j} = 0$ for $i \not =j$. For $i=j$ we get
\begin{eqnarray*}
\Lambda_{i}=\Lambda_{i,i} & = &
- \left( \begin{array}{ccc}
 v_{i,1} &  \ldots & v_{i,m_{i}} \\
 \end{array} \right ) \, {K_{i}}^{-1} \, {K_{i}} \, {K_{i}}^{-1} \,
\left( \begin{array}{c}
  v_{i,1}^{T} \\ \vdots \\  v_{i,m_{i}}^{T} \\
  \end{array}  \right) \\
& = &
 \left( \begin{array}{ccc}
 v_{i,1} &  \ldots & v_{i,m_{i}} \\
 \end{array} \right ) \, {(-K_{i})}^{-1} \,
\left( \begin{array}{c}
  v_{i,1}^{T} \\ \vdots \\  v_{i,m_{i}}^{T} \\
  \end{array}  \right),
\end{eqnarray*}
where $$ - K_{i} =
 \left( \begin{array}{c}
  v_{i,1}^{T} \\ \vdots \\  v_{i,m_{i}}^{T} \\
  \end{array}  \right) \, K_{n-1}(x_{i},x_{i}) \,
 \left( \begin{array}{ccc}
 v_{i,1} \ \ldots \ v_{i,m_{i}} \\
 \end{array} \right ),  $$
which was to be proved. \\
Now we show that this quadrature formula is exact for
matrix polynomials $F(x)$ and $G(x)$ satisfying
$$ \deg{F(x)}+\deg{G(x)} \leq 2n-1.$$
If $F(x)$ and $G(x)$ are matrix polynomials with degree $ \leq n-1$, then the
quadrature formula will be exact because the interpolation matrix polynomials
are exactly $F(x)$ and $G(x)$.
Let $F(x)$ be a matrix poynomial of degree $n+l$ and $G(x)$ a matrix
polynomial of degree $n-l-1$, where $l=0,1,\ldots,n-1$.
Then $G(x)$ satisfies the equation
$$ G(x) = \sum_{i=1}^{k} G(x_{i})W_{i}(x), $$
because $\deg{G(x)} \leq n-1$.
On the other hand $F(x)$ is approximated by a matrix polynomial of degree $n-1$ :
$$ F(x) \simeq \sum_{i=1}^{k} F(x_{i})W_{i}(x) = P(x). $$
This implies that $F(x)-P(x)$ is a matrix polynomial of degree $n+l$,
with zeros $x_{1},x_{2},\ldots,x_{k}$ and rootvectors
$v_{1,1},\ldots,v_{1,m_{1}},\ldots,v_{k,1},\ldots,v_{k,m_{k}}$. But then we
have  $$ F(x) - P(x) = R(x) \, P_{n}(x), $$
with $R(x)$ a matrix polynomial of degree $l$, (see Section 1).
So we get
\begin{eqnarray*}
 \int_{a}^{b} F(x) \, W(x) \, G(x)^{T}\, dx   & = &
 \int_{a}^{b} P(x) \, W(x) \, G(x)^{T}\, dx  \, + \,
 \int_{a}^{b} R(x)P_{n}(x) \, W(x) \, G(x)^{T}\, dx  \\
 & = &
\sum_{i=1}^{k} P(x_{i})\, \Lambda_{i} \, G(x_{i})^{T}  \, + \,
\int_{a}^{b} \left( \sum_{i=0}^{l} R_{i}x^{i} \right)\, P_{n}(x) \, W(x) \,
G(x)^{T}\, dx \\ & = &
\sum_{i=1}^{k} P(x_{i})\, \Lambda_{i} \, G(x_{i})^{T}  \, + \,
\sum_{i=0}^{l} R_{i} \int_{a}^{b} x^{i}\, P_{n}(x) \, W(x) \, G(x)^{T}\, dx
\\ & = &
\sum_{i=1}^{k} P(x_{i})\, \Lambda_{i} \, G(x_{i})^{T}  \, + \,
\sum_{i=0}^{l} R_{i} \int_{a}^{b} P_{n}(x) \, W(x) \, x^{i}G(x)^{T}\, dx
\end{eqnarray*}
The last integral will vanish because $ \deg{x^{i}G(x)} \leq l+n-l-1 = n-1
$ and $P_{n}(x)$ is orthogonal to every matrix polynomial of lower degree.
Further we have
\begin{eqnarray*}
P(x_{i})\, \Lambda_{i} & = &
P(x_{i})\,
\left( \begin{array}{ccc}
 v_{i,1} &  \ldots & v_{i,m_{i}} \\
 \end{array} \right ) \, {L_{i}}^{-1} \,
\left( \begin{array}{c}
  v_{i,1}^{T} \\ \vdots \\  v_{i,m_{i}}^{T} \\
  \end{array}  \right)  \\
& = & F(x_{i})\,
\left( \begin{array}{ccc}
 v_{i,1} &  \ldots & v_{i,m_{i}} \\
 \end{array} \right ) \, {L_{i}}^{-1} \,
\left( \begin{array}{c}
  v_{i,1}^{T} \\ \vdots \\  v_{i,m_{i}}^{T} \\
  \end{array}  \right),  \\
\end{eqnarray*}
so that
$$  \int_{a}^{b} F(x) \, W(x) \, G(x)^{T}\, dx    =
\sum_{i=1}^{k} F(x_{i})\, \Lambda_{i} \, G(x_{i})^{T}. $$
The case where $\deg{G(x)} \geq n$ is treated in a similar way, but we can
also take the transpose of the above mentioned equation.
\endproof\bigskip

\noindent
We showed that using the zeros and the rootvectors of the orthonormal
matrix polynomials on the interval $[a,b]$ with respect
to the weight matrix function $W(x)$ leads to a quadrature formula
with degree of precision $2n-1$.

\noindent
Note that $L_{i}$ is a symmetric and positive definite matrix
(see proof of Theorem \ref{eig:intlagr}) and so $\Lambda_{i}$ is a symmetric
and positive semidefinite matrix.

\noindent
\underline{Example.} \\
Consider the orthonormal matrix polynomials on the interval $[-1,1]$ with
respect to the weight matrix function
$$ W(x) = \left(
\begin{array}{cc} \frac{1}{\pi} {(1-x^{2})}^{-1/2}  &     0   \\
 0   & \frac{2}{\pi} {(1-x^{2})}^{1/2}    \\
\end{array} \right ). $$
We already mentioned that
$$ P_{2}(x) =
 \left( \begin{array}{cc}
 2 \sqrt{2} x^{2} - \sqrt{2}  & 0   \\
 0  & 4 x^{2}-1   \\
\end{array} \right ), $$
with zeros
$$ x_{1}=\frac{1}{\sqrt{2}}, \quad \quad x_{2}=\frac{-1}{\sqrt{2}}, \quad \quad
x_{3}=\frac{1}{2}, \quad \quad x_{4}=\frac{-1}{2} $$
and rootvectors
$$ v_{1}= \left( \begin{array}{c}  1  \\  0 \end{array} \right), \quad
   v_{2}= \left( \begin{array}{c}  4  \\  0 \end{array} \right), \quad
   v_{3}= \left( \begin{array}{c}  0  \\  3 \end{array} \right), \quad
   v_{4}= \left( \begin{array}{c}  0  \\ -2 \end{array} \right).$$
The quadrature coefficients can be computed as follows :
\begin{eqnarray*}
\Lambda_{1} & = &
 \left( \begin{array}{c}  1  \\  0 \end{array} \right) \,
 {\left( {
 \left( \begin{array}{cc} 1 & 0  \\  \end{array} \right) \,
\left( \begin{array}{cc}  1+1 & 0 \\ 0 & 1+2 \\ \end{array} \right) \,
 \left( \begin{array}{c}  1 \\  0  \\ \end{array} \right) }
 \right) }^{-1} \,
 \left( \begin{array}{cc} 1 & 0 \\  \end{array} \right)   \\
 & = &
 \left( \begin{array}{c}  1  \\  0 \end{array} \right) \, \frac{1}{2} \,
 \left( \begin{array}{cc} 1 & 0 \\  \end{array} \right)  \\
 & = &
 \left( \begin{array}{cc}
1/2  &  0 \\   0  &  0  \\
\end{array} \right).
\end{eqnarray*}
In a similar way we find
$$ \Lambda_{2}  =
 \left( \begin{array}{cc} 1/2  &  0 \\   0  &  0  \\
\end{array} \right) \quad \quad
 \Lambda_{3}  =
 \left( \begin{array}{cc} 0  &  0 \\   0  &  1/2  \\
\end{array} \right) \quad \quad
 \Lambda_{4}  =
 \left( \begin{array}{cc} 0  &  0 \\   0  &  1/2  \\
\end{array} \right). $$
Suppose we have the following matrix polynomials
 $$ F(x) = \left( \begin{array}{cc}
 x^{2}+1  & 6x   \\
 7x+1    & 5x^{2}-1  \\
\end{array} \right ) \quad \quad {\rm and } \quad \quad
G(x) = \left( \begin{array}{cc}
 2x+5  & 6x   \\
   7   & 4x-3  \\
\end{array} \right ), $$
then
\begin{eqnarray*}
\int_{a}^{b} F(x) \, W(x)\, G(x)^{T} \, dx  & = &
 \left( \begin{array}{cc}
 \frac{3}{2}  & 3\sqrt{2}   \\
 \frac{7}{\sqrt{2}}+1  & \frac{3}{2}   \\
\end{array} \right ) \, \Lambda_{1}\,
 \left( \begin{array}{cc}
 5+\sqrt{2}  & 7   \\
  3\sqrt{2}   & 2 \sqrt{2}-3  \\
\end{array} \right ) \\
& + &  \left( \begin{array}{cc}
 \frac{3}{2}  & -3\sqrt{2}   \\
\frac{-7}{\sqrt{2}}+1  & \frac{3}{2}   \\
\end{array} \right ) \, \Lambda_{2} \,
 \left( \begin{array}{cc}
 5-\sqrt{2}  & 7   \\
 -3\sqrt{2}   & -2 \sqrt{2}-3  \\
\end{array} \right )  \\
& + & \left( \begin{array}{cc}
\frac{5}{4}  & 3   \\
\frac{9}{2} & \frac{1}{4}  \\
\end{array} \right ) \, \Lambda_{3} \,
 \left( \begin{array}{cc}
 6  & 7   \\
  3 & -1  \\
\end{array} \right ) \, + \,
 \left( \begin{array}{cc}
\frac{5}{4}  & -3   \\
\frac{-5}{2} & \frac{1}{4}  \\
\end{array} \right ) \, \Lambda_{4} \,
 \left( \begin{array}{cc}
 4  & 7   \\
 -3 & -5  \\
\end{array} \right ) \\
& = & \left( \begin{array}{cc}
33/2  &   33/2   \\
   12   & 25/4 \\
\end{array} \right ).\\
\end{eqnarray*}

\subsection{Convergence}
The Gaussian quadrature formula converges to the exact value of the matrix
integral of the functions $F(x)$ and $G(x)$, even without imposing severe
conditions on those functions.

\begin{theorem}
\label{eig:conv}
Let $F(x)$ and $G(x)$ be continuous matrix functions on the finite interval
$[a,b]$, then we have
$$ \int_{a}^{b} F(x) \, W(x) \, G(x)^{T}\, dx = \lim_{n \rightarrow \infty}
\sum_{i=1}^{k} F(x_{i}^{(n)})\, \Lambda_{i}^{(n)} \, {G(x_{i}^{(n)})}^{T}, $$
with $x_{i}^{(n)}$ the zeros of the orthonormal matrix polynomial
$P_{n}(x)$ on the interval $[a,b]$ with respect to the weight matrix
function $W(x)$ and with $\Lambda_{i}^{(n)}$ the corresponding
quadrature coefficients.
\end{theorem}

\noindent
In the proof of this theorem we will use the following proposition of
matrix norms:
\begin{proposition}
If $ A_{1},A_{2},\ldots,A_{k}$ are symmetric positive semidefinite $p
\times p$-matrices, then they satisfy
$$\sum_{i=1}^{k} {\|A_{i}\|}_{2}  \leq p \,
{\left \| \sum_{i=1}^{k} A_{i} \right \| }_{2}, $$
with $$ {\|A\|}_{2} = \sqrt{ \rho (A^{*}A) } = \rho (A)$$
and $\rho (A) $ is the spectral radius of $A$.
\end{proposition}

\noindent
\underline{Proof}\\
Suppose $ \alpha_{1}^{(i)} \geq \ldots \geq \alpha_{p}^{(i)} \geq 0 $ are the
eigenvalues of $A_{i}$ and
$ \gamma_{1} \geq \ldots \geq \gamma_{p} \geq 0 $ are those of
$\sum_{i=1}^{k} A_{i}$. Then we have
\begin{eqnarray*}
 \sum_{i=1}^{k} {\|A_{i}\|}_{2}
& = & \sum_{i=1}^{k} \alpha_{1}^{(i)} \, \leq \,
\sum_{i=1}^{k} \sum_{j=1}^{p} \alpha_{j}^{(i)}
\, = \, \sum_{i=1}^{k} tr(A_{i}) \, = \,   tr( \sum_{i=1}^{k} A_{i} )  \\
& = & \sum_{j=1}^{p} \gamma_{j} \, \leq \, p \, \gamma_{1}
\, = \, p \, {\left \| \sum_{i=1}^{k} A_{i} \right \|}_{2}.\ \ \rule{0.5em}{0.5em}
\end{eqnarray*}
\bigskip

\noindent
If the interval $[a,b]$ is finite we have another useful inequality, namely
$$ \| \int_{a}^{b} F(x) \, d \delta (x) \| \, \leq \,
 \int_{a}^{b} \|F(x)\| \, d \delta (x), $$
(see \cite[p. 78]{ru:1973}).

\noindent
\underline{Proof of Theorem \ref{eig:conv}} \\
We will distinguish 3 cases :
\begin{itemize}
\item[(1)] $ G(x) = I $.  \\
Since the matrix function $F(x)$ is continuous, every $F_{i,j}(x)$
will be continuous.
This means that every element can be approximated arbitrarily close by a
polynomial and thus
$$ \| F(x) - Q_{m}(x) \|_{2} \leq \varepsilon. $$
So we have
\begin{eqnarray*}
\lefteqn{
\| \int_{a}^{b} F(x) \, W(x) \, dx -
\sum_{i=1}^{k} F(x_{i}^{(n)})\, \Lambda_{i}^{(n)} \|_{2} } \\
& \leq &
\|\int_{a}^{b} F(x) \, W(x) \, dx - \int_{a}^{b} Q_{m}(x) \, W(x) \, dx \|_{2}
\\
& + &
\| \int_{a}^{b} Q_{m}(x) \, W(x) \, dx -
\sum_{i=1}^{k} Q_{m}(x_{i}^{(n)})\, \Lambda_{i}^{(n)} \|_{2}  \\
& + &
\| \sum_{i=1}^{k} Q_{m}(x_{i}^{(n)})\, \Lambda_{i}^{(n)} -
\sum_{i=1}^{k} F(x_{i}^{(n)})\, \Lambda_{i}^{(n)} \|_{2}.
\end{eqnarray*}
But the first term of this sum satisfies
\begin{eqnarray*}
\| \int_{a}^{b} F(x) \, W(x) \, dx - \int_{a}^{b} Q_{m}(x) \, W(x) \, dx
\|_{2} & \leq &
\int_{a}^{b}  {\|F(x)-Q_{m}(x)\|_{2}}\, {\|W(x)\|_{2}} \, dx \\
& \leq &
\varepsilon \int_{a}^{b}   {\|W(x)\|_{2}} \, dx. \\
\end{eqnarray*}
The second term
$$ \| \int_{a}^{b} Q_{m}(x) \, W(x) \, dx -
\sum_{i=1}^{k} Q_{m}(x_{i}^{(n)})\, \Lambda_{i}^{(n)} \|_{2} = 0,  $$
for $m \leq 2n-1$ and the third term can be bounded as follows :
\begin{eqnarray*}
\| \sum_{i=1}^{k} Q_{m}(x_{i}^{(n)})\, \Lambda_{i}^{(n)} -
\sum_{i=1}^{k} F(x_{i}^{(n)})\, \Lambda_{i}^{(n)} \|_{2}
& \leq &
 \sum_{i=1}^{k} {\|Q_{m}(x_{i}^{(n)})-
 F(x_{i}^{(n)})\|_{2}}\, {\|\Lambda_{i}^{(n)}\|_{2}}   \\
& \leq &
\varepsilon  \sum_{i=1}^{k} {\|\Lambda_{i}^{(n)}\|_{2}}.
\end{eqnarray*}
So we have
$$ \| \int_{a}^{b} F(x) \, W(x) \, dx -
\sum_{i=1}^{k} F(x_{i}^{(n)})\, \Lambda_{i}^{(n)} \|_{2} \leq
\varepsilon \,
 \left( \int_{a}^{b} {\|W(x)\|_{2}}\, dx \, + \,
 \sum_{i=1}^{k}  {\| \Lambda_{i}^{(n)} \|_{2}} \right). $$
But from the previous proposition we get
$$ \sum_{i=1}^{k}  {\| \Lambda_{i}^{(n)} \|_{2}} \leq
p \,  \| \sum_{i=1}^{k} \Lambda_{i}^{(n)} \|_{2} =
p \,  \| \int_{a}^{b} W(x)\, dx  \|_{2} \leq
p \,  \int_{a}^{b} \|W(x)\|_{2} \, dx  $$
and thus we have
$$ \| \int_{a}^{b} F(x) \, W(x) \, dx -
\sum_{i=1}^{k} F(x_{i}^{(n)})\, \Lambda_{i}^{(n)} \|_{2} \leq
\varepsilon \, (1+p) \, \int_{a}^{b} \|W(x)\|_{2} \, dx.  $$
But this means that
$$ \int_{a}^{b} F(x) \, W(x) \, dx = \lim_{n \rightarrow \infty}
\sum_{i=1}^{k} F(x_{i}^{(n)})\, \Lambda_{i}^{(n)}. $$

\item[(2)] $G(x)$ is a matrix polynomial.  \\
Suppose $G(x) = G_{l}x^{l}+G_{l-1}x^{l-1}+\ldots+G_{1}x+G_{0}$, then we have
\begin{eqnarray*}
\int_{a}^{b} F(x) \, W(x) \, G(x)^{T}\, dx & = &
\sum_{j=0}^{l}\, \left( \int_{a}^{b} F(x) \, W(x) \, x^{j}\, dx \right) \,
G_{j}^{T} \\
& = &
\sum_{j=0}^{l}\, \left( \int_{a}^{b} x^{j}F(x) \, W(x) \,  dx \right) \,
G_{j}^{T} \\
& = &
\sum_{j=0}^{l}\, \lim_{n \rightarrow \infty}
\sum_{i=1}^{k} {(x_{i}^{(n)})}^{j}\, F(x_{i}^{(n)})\, \Lambda_{i}^{(n)}
 G_{j}^{T}  \\
& = &
\lim_{n \rightarrow \infty}
\sum_{i=1}^{k}  F(x_{i}^{(n)})\, \Lambda_{i}^{(n)} \,
\sum_{j=0}^{l}\,{(x_{i}^{(n)})}^{j}\, G_{j}^{T}  \\
& = &
\lim_{n \rightarrow \infty}
\sum_{i=1}^{k}  F(x_{i}^{(n)})\, \Lambda_{i}^{(n)} \, G(x_{i}^{(n)})^{T}.
\end{eqnarray*}
Which implies that the theorem holds in the case that
$G(x)$ is a matrix polynomial.
\item[(3)] $G(x)$ is a continuous matrix function on the interval $[a,b]$.\\
In this case, both functions $F(x)$ and $G(x)$ can be approximated arbitrarily
closely by a matrix polynomial :
$$ \| F(x) - Q_{m}(x) \|_{2} \leq \varepsilon  \quad {\rm and } \quad
 \| G(x) - P_{l}(x) \|_{2} \leq \varepsilon. $$
This gives
\begin{eqnarray*}
\lefteqn{ \| \int_{a}^{b} F(x)\, W(x)\, G(x)^{T}\, dx  -
\sum_{i=1}^{k} F(x_{i}^{(n)})\, \Lambda_{i}^{(n)}\, {G(x_{i}^{(n)})}^{T} \|_{2}
} \\ & \leq &
\| \int_{a}^{b} F(x)\, W(x)\, G(x)^{T}\, dx -
\int_{a}^{b} Q_{m}(x)\, W(x)\, G(x)^{T}\, dx  \\
& & - \,
\int_{a}^{b} F(x)\, W(x)\, P_{l}(x)^{T}\, dx +
\int_{a}^{b} Q_{m}(x)\, W(x)\, P_{l}(x)^{T}\, dx \|_{2}  \\
& + &
\| \int_{a}^{b} Q_{m}(x)\, W(x)\, G(x)^{T}\, dx -
\sum_{i=1}^{k} Q_{m}(x_{i}^{(n)})\, \Lambda_{i}^{(n)}\, {G(x_{i}^{(n)})}^{T}
\|_{2} \\
& + &
\| \int_{a}^{b} F(x)\, W(x)\, P_{l}(x)^{T}\, dx -
\sum_{i=1}^{k} F(x_{i}^{(n)})\, \Lambda_{i}^{(n)}\, {P_{l}(x_{i}^{(n)})}^{T}
\|_{2} \\
& + &
\| \int_{a}^{b} Q_{m}(x)\, W(x)\, P_{l}(x)^{T}\, dx -
\sum_{i=1}^{k} Q_{m}(x_{i}^{(n)})\, \Lambda_{i}^{(n)}\, {P_{l}(x_{i}^{(n)})}^{T}
\|_{2} \\
& + &
\| \sum_{i=1}^{k} Q_{m}(x_{i}^{(n)})\, \Lambda_{i}^{(n)} {G(x_{i}^{(n)})}^{T} +
\sum_{i=1}^{k} F(x_{i}^{(n)})\, \Lambda_{i}^{(n)}\, {P_{l}(x_{i}^{(n)})}^{T} \\
& & - \,
\sum_{i=1}^{k} Q_{m}(x_{i}^{(n)})\, \Lambda_{i}^{(n)}\, {P_{l}(x_{i}^{(n)})}^{T}
- \sum_{i=1}^{k} F(x_{i}^{(n)})\, \Lambda_{i}^{(n)}\, {G(x_{i}^{(n)})}^{T}
\|_{2}. \end{eqnarray*}
From the previous case we get that the second and the third term are arbitrarily
small ($< \varepsilon$) for $n$ sufficiently large. The fourth term is equal to
$0$ if $m+l \leq 2n-1$ and the other two terms can be bounded as follows,
\begin{eqnarray*}
\lefteqn{
\| \int_{a}^{b} [F(x)-Q_{m}(x)]\, W(x)\, {[G(x)-P_{l}(x)]}^{T}\, dx \|_{2}} \\
& \leq &  \int_{a}^{b}\, \|F(x)-Q_{m}(x)\|_{2}\, \|W(x)\|_{2}\,
\|G(x)-P_{l}(x)\|_{2}\, dx \\
& \leq & {\varepsilon}^{2}\,  \int_{a}^{b}\, \|W(x)\|_{2}\, dx \\
\end{eqnarray*}
and
\begin{eqnarray*}
\lefteqn{
\| \sum_{i=1}^{k} [F(x_{i}^{(n)})- Q_{m}(x_{i}^{(n)})]\, \Lambda_{i}^{(n)}\,
{[G(x_{i}^{(n)})-P_{l}(x_{i}^{(n)})]}^{T}\|_{2}}\\
& \leq &  \sum_{i=1}^{k}
\|F(x_{i}^{(n)})- Q_{m}(x_{i}^{(n)})\|_{2}\, \|\Lambda_{i}^{(n)}\|_{2}\,
\|G(x_{i}^{(n)})-P_{l}(x_{i}^{(n)})\|_{2} \\
& \leq & {\varepsilon}^{2}  \sum_{i=1}^{k}
\|\Lambda_{i}^{(n)}\|_{2}.
\end{eqnarray*}
Hence for $n$ large enough
\begin{eqnarray*}
\lefteqn{ \|\int_{a}^{b} F(x)\, W(x)\, G(x)^{T}\, dx  -
\sum_{i=1}^{k} F(x_{i}^{(n)})\, \Lambda_{i}^{(n)}\, G(x_{i}^{(n)})^{T} \|_{2} }
\\
& \leq & 2 \varepsilon + {\varepsilon}^{2} \,
 \left( \int_{a}^{b} {\|W(x)\|_{2}}\, dx \, + \,
 \sum_{i=1}^{k}  {\| \Lambda_{i}^{(n)} \|_{2}} \right) \\
 & \leq & 2 \varepsilon + {\varepsilon}^{2} \, (1+p) \,
\int_{a}^{b} {\|W(x)\|_{2}}\, dx,
\end{eqnarray*}
which proves the convergence of the Gaussian quadrature formula.
\endproof
\end{itemize}

\end{document}